\documentclass[a4paper,10pt]{article}

\usepackage{latexsym}
\usepackage{mathrsfs}
\usepackage{amsfonts,amsmath,amssymb}
\usepackage{indentfirst}
\usepackage{subeqnarray}
\usepackage{verbatim}
\usepackage[pdftex]{graphicx}
\usepackage{epstopdf}
\usepackage{stmaryrd}
\usepackage{multirow}
\usepackage{diagbox}
\usepackage{subcaption}
\usepackage{float}
\usepackage{booktabs}

\newcommand{\bx}{\mbox{\boldmath{$x$}}}

\newcommand{\bb}{\mbox{\boldmath{$b$}}}

\newcommand{\bn}{\mbox{\boldmath{$n$}}}

\newcommand{\bW}{\mbox{\boldmath{$W$}}}

\newcommand{\bs}{\mbox{\boldmath{$s$}}}

\newcommand{\dx}{\mathrm{d}x}
\newcommand{\ds}{\mathrm{d}s}
\newcommand{\dt}{\mathrm{d}t}

\newtheorem{theorem}{Theorem}[section]

\newtheorem{example}[theorem]{Example}

\newtheorem{remark}[theorem]{Remark}
\numberwithin{equation}{section}

\usepackage[pdftex,unicode,colorlinks,linkcolor=blue]{hyperref}

\begin{document}

\begin{center}
\Large\bf Local Randomized Neural Networks with Discontinuous Galerkin Methods for KdV-type and Burgers Equations
\end{center}

\begin{center}
Jingbo Sun\footnote{School of Mathematics and Statistics, Xi'an Jiaotong University, Xi'an, Shaanxi 710049, P.R. China. E-mail: {\tt jingbosun@stu.xjtu.edu.cn}.},
\quad 
Fei Wang\footnote{School of Mathematics and Statistics, Xi’an Jiaotong University, Xi’an, Shaanxi 710049, China. The work of this author was partially supported by the National Natural Science Foundation of China (Grant No. 12171383). Email: {\tt feiwang.xjtu@xjtu.edu.cn}.}
\end{center}

\medskip
\begin{quote}
  {\bf Abstract.} The Local Randomized Neural Networks with Discontinuous Galerkin (LRNN-DG) methods, introduced in \cite{Sun2022lrnndg}, were originally designed for solving linear partial differential equations. In this paper, we extend the LRNN-DG methods to solve nonlinear PDEs, specifically the Korteweg-de Vries (KdV) equation and the Burgers equation, utilizing a space-time approach. Additionally, we introduce adaptive domain decomposition and a characteristic direction approach to enhance the efficiency of the proposed methods. Numerical experiments demonstrate that the proposed methods achieve high accuracy with fewer degrees of freedom, additionally, adaptive domain decomposition and a characteristic direction approach significantly improve computational efficiency.

\end{quote}

\medskip

\section{Introduction}

This research explores the application of Local Randomized Neural Networks with Discontinuous Galerkin methods for solving two important nonlinear wave equations: the Korteweg-de Vries (KdV) equation and the Burgers equation. The KdV equation, fundamental for studying wave behavior on shallow water surfaces, is well-known for its nonlinear nature and third-order derivative terms (\cite{Korteweg1989KdV}). Similarly, the Burgers equation, originally introduced in turbulence modeling by Burgers (\cite{Burgers2013bur}), has widespread applications in physics and engineering. Several numerical methods, such as finite difference, finite element, discontinuous Galerkin, and spectral methods, have been developed to solve these equations (\cite{Zabusky1965KdVff1, Vliengenthart1971ff2, Cheng2015mfekdvbur, Ma1986spectralkdv, Maday1988spectralkdv, Yan2002dgkdv, Xu2007dgkdv, Mittal1993nsBur, Tan2022dgkdv, Ozis2003fembur, Dogan2004gfemBur, Caldwell1982spebur}).

Recently, neural network-based methods such as the Deep Ritz Method (\cite{Ee2018deepRitz}), Deep Galerkin Method (\cite{Sirignano2018DGM}), Physics-Informed Neural Networks (\cite{Raissi2019PINN}), Weak Adversarial Networks (\cite{Zang2020adversarialnns}), and Deep Nitsche Method (\cite{Liao2019deepRitzboundary}) have attracted significant interest from researchers. These methods leverage the powerful approximation capabilities of neural networks, which are theoretically supported by studies such as \cite{Cybenko1989appro, Chen1995appro, Hornik1991appro, Mhaskar1995appro, Barron1993appro, Ma2022appro, Lu2017appro}. However, a major limitation of neural network-based approaches is their reliance on optimization solvers for training, which can lead to difficulties in achieving high accuracy and efficiency, largely due to the challenges of solving nonconvex optimization problems. Consequently, traditional numerical methods often outperform neural networks in terms of both accuracy and computational efficiency.

To address the challenges of nonconvex optimization, randomized neural networks (RNNs) have been proposed (\cite{Igelnik1995RNN1, Igelnik1999RNN2, Pao1992RNN3, Pao1994RNN4}). RNNs differ from conventional neural networks in their training process. In RNNs, the weights between hidden layers are randomly assigned and kept fixed during training, while the output layer parameters are determined using a least-squares approach. Studies on RNNs' approximation capabilities, including \cite{Igelnik1995RNN1, Liu2014ELMfeasible, Neufeld2023rnn, Huang2006ELMtheorandapp}, have shown that RNNs can achieve comparable approximation errors to standard neural networks, provided that activation functions and parameter initialization strategies are carefully chosen.

Building on these ideas, Dong and Li introduced the Local Extreme Learning Machine (ELM), which integrates ELM with non-overlapping domain decomposition techniques (\cite{Dong2021locELM}). Further exploration by Dong and Wang examined the influence of initialization on this method in \cite{Dong2021locELMW}. Additionally, Shang et al. proposed the RNN-Petrov-Galerkin (RNN-PG) method for solving both linear and nonlinear PDEs, integrating RNNs with the Petrov-Galerkin formulation (\cite{Shang2022DeepPetrov,Shang2023RNNPG}). Sun et al. introduced the Local Randomized Neural Network with Discontinuous Galerkin (LRNN-DG) methods for solving linear PDEs in \cite{Sun2022lrnndg, Sun2024lrnndgdvwe}. These methods combine domain decomposition techniques and utilize distinct local RNNs to approximate numerical solutions in each sub-domain, with the DG method employed to couple the local solutions. Numerical experiments have demonstrated that LRNN-DG methods can effectively solve time-dependent problems with greater accuracy and fewer degrees of freedom compared to traditional methods like the discontinuous Galerkin approach. As space-time methods, LRNN-DG obtains numerical solutions via a least-squares solver, reducing error accumulation across time iterations. Other works on RNN-based methods can be found in \cite{Chen2022rfm,Zhang2024Transferable,Dang2023lrnnwg,Li2023lrnnip,Dang2024adaptive} and the references therein.

While LRNN-DG methods have been applied to solve linear problems such as the Poisson equation, heat equation, and diffusive-viscous wave equation (\cite{Sun2022lrnndg, Sun2024lrnndgdvwe}), they have not yet been extended to nonlinear PDEs. This paper develops a space-time LRNN-DG method to solve nonlinear KdV and Burgers equations, incorporating suitable mesh generation strategies to further enhance the performance of the neural network.

The paper is organized as follows: Section \ref{sec:RNN} provides an overview of the architecture and training approach of RNNs, along with the notations for the DG formulation. Section \ref{dgforkdv} introduces the LRNN-DG method for solving the KdV equation, while Section \ref{dgforburgers} details its application to Burgers equations. Section \ref{sec:adap} presents the adaptive and characteristic meshes used in the experiments. Section \ref{sec:ex} showcases the numerical results of applying LRNN-DG methods on different mesh types. Finally, Section \ref{summary} concludes the paper with remarks and future research directions.

%%%%%%%%%%%%%%%%%%%%%%%%%%%%%
\section{Network Structure and Notation}
\label{sec:RNN}

In this section, we introduce the concept of randomized neural networks used in this study, along with relevant notation related to the Discontinuous Galerkin (DG) method.

\subsection{Randomized Neural Networks}

Let $I$ represent the time interval, and let $\Omega\subset\mathbb{R}^d$ denote the spatial domain of interest. The input vector is $\bs\in \Sigma:= I\times \Omega$. The fully connected neural network structure is defined as follows:
\begin{subequations}\label{fnn_structure}
\begin{align}
N^{(1)}(\bs) &=\rho(\bW^{(1)} \bs + \bb^{(1)}),\\
N^{(i)}(\bs) &=\rho(\bW^{(i)} N^{(i-1)} + \bb^{(i)}),\quad i= 2,\cdots, L,\\
\mathcal{U}(\bs) &=\bW^{(L+1)} N^{(L)},
\end{align}
\end{subequations}
where $N^{(i)}$ denotes the $i$-th hidden layer with the weight matrix $\bW^{(i)}$ and bias vector $\bb^{(i)}$. The activation function $\rho$ can be either Tanh or ReLU. $L$ represents the number of hidden layers (depth), and $\bW^{(L+1)}$ corresponds to the weight matrix of the output layer, with the bias term in the output layer omitted. The set of all functions representable by this network is denoted as:
\begin{equation*}
\mathcal{M}(\theta,L,\Sigma)=\{\mathcal{U}(\bs)= \bW^{(L+1)} (N^{(L)} \circ \cdots \circ N^{(1)}(\bs)),\ \bs\in \Sigma\},
\end{equation*}
where $\theta = \{\bW^{(L+1)},(\bW^{(l)} , \bb^{(l)})^{L}_{l=1}\}$ and $\circ$ denotes function composition.

\begin{figure}[thpb] 
  \centering  
  \includegraphics[width=0.8\textwidth]{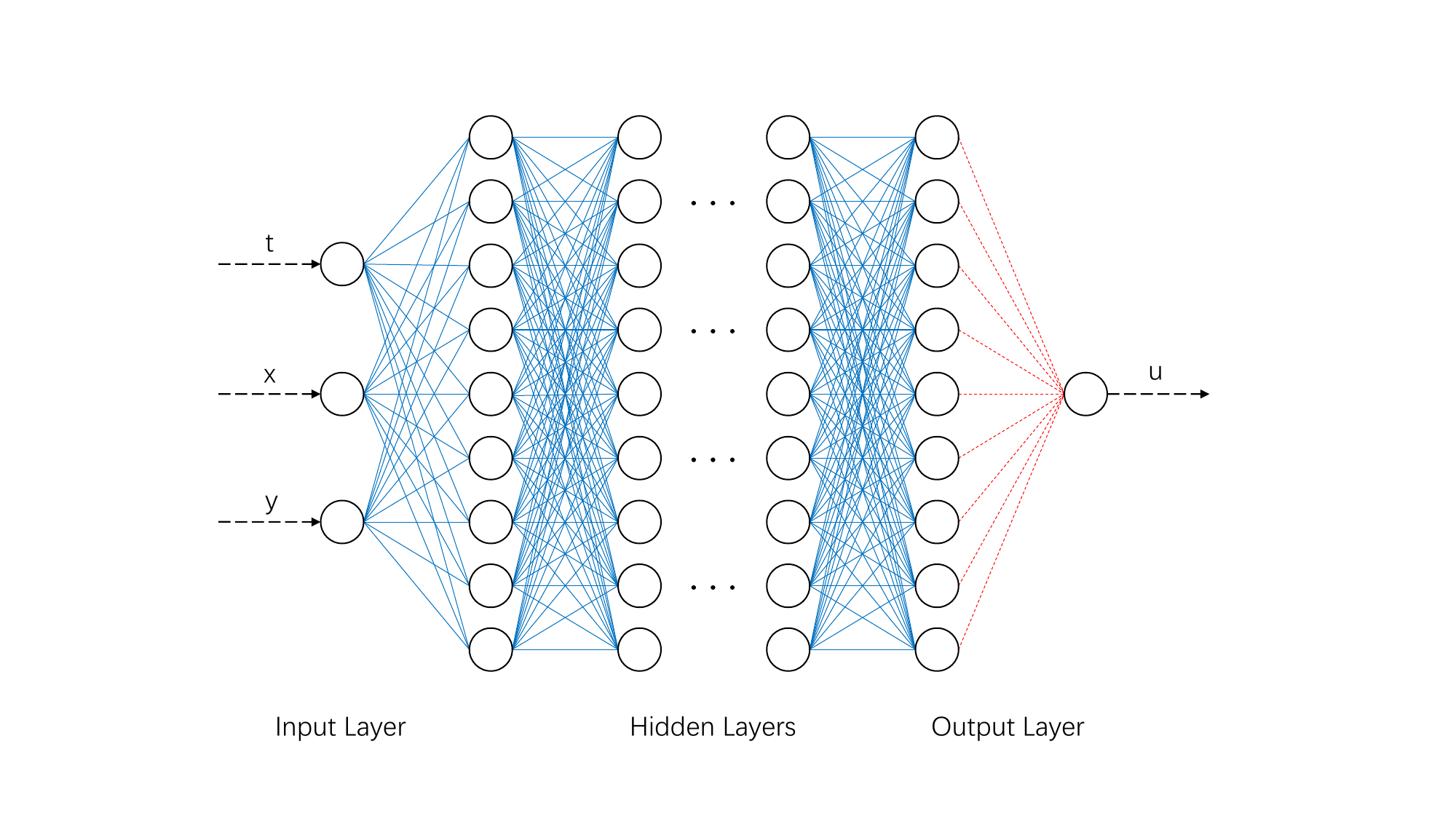}  
  \caption{The structure of a randomized neural network}
  \label{nn_structure}
\end{figure}

Training all parameters $\theta$ in the fully connected neural network involves solving a nonlinear, nonconvex optimization problem, which is computationally expensive and prone to local minima. In randomized neural networks, except for the parameters between the last hidden layer and the output layer (computed using least-squares methods), all other parameters are randomly initialized and kept fixed throughout the training process. To illustrate the concept of randomized neural networks more clearly, see Figure \ref{nn_structure}. In this figure, the input consists of $t$, $x$, and $y$, while the output is the solution $u$. The parameters along the blue solid lines are randomly assigned and fixed, whereas the parameters along the red dashed lines are determined using the least-squares method.

\vspace{-2mm}

Due to the absence of a bias term in the output layer, the function represented by the neural network can be viewed as a linear combination of nonlinear basis functions. This leads to an alternative representation of $\mathcal{M}(\theta,L,\Omega)$, given by:
\begin{equation}\label{rnn_space}
\mathcal{M}(\sigma)=\left\{\mathcal{U}(\alpha, \theta, t, \bx) = \sum^M_{j=1} \alpha^\sigma_j \phi^\sigma (\theta_j, t, \bx) : 
\begin{pmatrix} t\\ \bx \end{pmatrix}\in \sigma\right\},
\end{equation}
where $\sigma\subset I \times\Omega$, and the nonlinear functions $\phi^\sigma (\theta_j, t, \bx)$ with parameters $\theta_j$ ($j=1,\cdots,M$) represent outputs of the last hidden layer. The parameters $\theta_j$ are randomly sampled from probability distributions, such as the uniform or Gaussian distribution, and are fixed. The set $\{\alpha^\sigma_j\}_{j=1}^M$ denotes the weights of the output layer, which are computed using least-squares methods, and $M$ represents the number of neurons in the last hidden layer. For simplicity, we write $\phi^\sigma_j(t,\bx)$ instead of $\phi^\sigma (\theta_j, t,\bx)$ in the subsequent sections.

\begin{remark}

Deep neural networks are known for their exceptional approximation capabilities. However, existing methods that rely on training neural networks often fall short in terms of accuracy when compared to traditional numerical methods. This shortcoming is primarily due to the significant optimization errors associated with training the network weights. Randomized neural networks mitigate this issue by reducing optimization complexity, albeit with a slight trade-off in approximation power. This results in a favorable balance between optimization error and approximation error, enabling RNNs to achieve considerably improved accuracy with lower computational cost.

\end{remark}

\subsection{Notation for Space-Time DG Formulation}

To address complex problems that may require multiple neural networks for effective solutions, different randomized neural networks can be employed to approximate solutions in various subdomains. These networks are integrated using the Discontinuous Galerkin (DG) scheme.

We first introduce some notation that will be used to elucidate the process of constructing a space-time Discontinuous Galerkin formulation. The spatial domain $\Omega$ is decomposed into a mesh $\{ \mathcal{T}_h \}$, where $h = \max_{K\in {\cal T}_h}\{{\rm diam}(K)\}$ and the number of elements in $\mathcal{T}_h$ is denoted as $N_s$. In the case of a one-dimensional spatial domain $\Omega = (x_0, x_{N_s})$, $\mathcal{T}_h$ consists of elements $K_i =(x_{i-1},x_i)$, for $i = 1, 2, \cdots, N_s$. Similarly, the temporal interval $I$ is divided into $N_t$ sub-intervals $\mathcal{D}_\tau =\{I_i =(t_{i-1},t_i), 0=t_0< t_1<\cdots<t_{N_t}=T\}$, where $\tau = \max\limits_{I_i\in {\cal D}_\tau}\{{\rm diam}(I_i)\}$, $t_0$ represents the initial time and $T$ denotes the final time. Furthermore, $\mathcal{E}_h$ denotes the union of all spatial mesh faces ( or edges), $\mathcal{E}_h^i$ represents the set of interior faces (edges), and $\mathcal{E}^\partial_h = \mathcal{E}_h\backslash {\cal E}_h^i$ signifies the set of boundary (faces) edges. The set of time nodes is denoted as ${\cal P}_\tau=\{t_i,i=0,\cdots,{N_t}\}$, and $\mathcal{P}_\tau^{i} = \mathcal{P}_\tau \backslash \{t_0,t_{N_t}\}$ refers to the set of all interior points. Consequently, we define the space-time decomposition $\mathcal{M}_{\tau h} = \mathcal{D}_\tau\times\mathcal{T}_h$ over the domain $\Sigma = I\times \Omega$, comprising a total of $N_e=N_t N_s$  sub-elements.

For adjacent elements $\sigma^+_h = I_i\times K^+$ and $\sigma^-_h= I_i\times K^-$ which share a common spatial face $f_h$, let $\boldsymbol{n}^{\pm}=\boldsymbol{n}|_{\partial K^{\pm}}$ be the unit outward normal vectors
on $\partial K^{\pm}$. For a scalar function $v$ and a
vector function $\boldsymbol{q}$, we define $v^{\pm}=v|_{\partial \sigma_h^{\pm}}$ and $\boldsymbol{q}^\pm=\boldsymbol{q}|_{\partial \sigma_h^{\pm}}$. The averages $\{ \cdot \}$ and jumps $\llbracket \cdot \rrbracket$, $[ \cdot ]$ on $f_h \in (\mathcal{D}_\tau\times\mathcal{E}_h^{i})$ are defined as:
\begin{subequations}
\begin{align*}
\{ v \} = \frac{1}{2} (v^+ + v^-),& \quad \llbracket v\rrbracket = v^+ \boldsymbol{n}^+ + v^- \boldsymbol{n}^-, \\
\{ \boldsymbol{q} \} = \frac{1}{2}(\boldsymbol{q}^+ +\boldsymbol{q}^-),& \quad
[\boldsymbol{q}] = \boldsymbol{q}^+ \cdot \boldsymbol{n}^+ +\boldsymbol{q}^- \cdot \boldsymbol{n}^-.
\end{align*}
If $f_h \in (\mathcal{D}_\tau\times \mathcal{E}_h^{\partial})$, we set
\begin{subequations}
\begin{align*}
\llbracket v\rrbracket = v \boldsymbol{n},\quad
\{ \boldsymbol{q} \} = \boldsymbol{q},
\end{align*}
\end{subequations}
where $\boldsymbol{n}$ is the unit outward normal vector on $\partial\Omega$.

Consider two adjacent elements $\sigma^+_\tau = I_{i+1}\times K$ and $\sigma^-_\tau= I_i\times K$ sharing a common temporal face $f_\tau \in\{t_i\}\times \mathcal{T}_h$. Here, we define $w(t^{\pm}_i,\bx)=w(t_i,\bx)|_{\partial \sigma_\tau^{\pm}}$ for a scalar function $w$ and establish the averages $\{ \cdot \}$ and jumps $[ \cdot ]$ on $f_\tau \in (\mathcal{P}_\tau^{i}\times \mathcal{T}_h)$ as follows:
\begin{align*}
\{ w(t_i,\bx) \} = \frac{1}{2} \left(w(t_i^+,\bx) + w(t_i^-,\bx)\right),& \quad [w(t_i,\bx)] =w(t_i^+,\bx) - w(t_i^-,\bx).
\end{align*}
\end{subequations}
If $f_\tau \in (\mathcal{P}_\tau^{\partial}\times \mathcal{T}_h)$, we set
\begin{subequations}
\begin{align*}
[w(t_0,\bx)] = - w(t_0,\bx),\quad [w(t_{N_t},\bx)] =  w(t_{N_t},\bx),\quad
\{ w(t,\bx) \} = w(t,\bx).
\end{align*}
\end{subequations}

To enable communication between sub-networks on adjacent elements, another approach can be used. Collocation points are introduced on interior edges, boundary edges, and initial edges to enforce $C^k$ continuity, boundary conditions, and initial conditions, respectively. The selection of these collocation points is denoted as follows: $N^b_h$ points $P^b_h= \{(t^b_h,\bx^b_h) \in \mathcal{D}_\tau\times\mathcal{E}^\partial_h \}$ on boundary edges, (we choose $N^-_h$ points $P^-_h= \{(t^-_h,x^-_h) \in \mathcal{D}_\tau\times\{ x_0\} \}$ and $N^+_h$ points $P^+_h= \{(t^+_h,x^+_h) \in \mathcal{D}_\tau\times\{ x_{N_s}\} \}$ for one-dimensional spatial domains), $N^0_\tau$ points $P^0_\tau= \{(t^0_\tau,\bx^0_\tau) \in t_0 \times\mathcal{T}_h \}$ on initial edges, $N^i_h$ points $P^i_h= \{(t^i_h,\bx^i_h) \in \mathcal{D}_\tau\times\mathcal{E}^i_h \}$ on spatial interior edges, and $N^i_\tau$ points $P^i_\tau= \{(t^i_\tau,\bx^i_\tau) \in \mathcal{P}^i_\tau\times\mathcal{T}_h \}$ on temporal interior edges. Depending on the scenario, appropriate conditions are applied to the numerical solution at these points, introducing additional systems of equations to address varying regularity conditions.

Finally, considering the space-time domain decomposition and \eqref{rnn_space}, we define the local randomized neural network function spaces as follows:
\begin{align}
V_{h\tau}&=\{v_{h\tau} \in L^2(\Sigma): \;v_{h\tau} |_\sigma \in \mathcal{M}_{RNN}(\sigma)\quad \forall\,\sigma\in\mathcal{D}_\tau\times\mathcal{T}_h \},\label{scalspace}\\
\boldsymbol{Q}_{h\tau}&=\{\boldsymbol{q}_{h\tau} \in [L^2(\Sigma)]^d: \;\boldsymbol{q}_{h\tau} |_\sigma \in \left[\mathcal{M}_{RNN}(\sigma)\right]^d\quad \forall\,\sigma\in\mathcal{D}_\tau\times\mathcal{T}_h \}\label{vecspace}.
\end{align}

\section{LRNN-DG Methods for KdV Equations}
\label{dgforkdv}

The KdV equation is a significant time-dependent nonlinear problem in physics. This section presents the development of two LRNN-DG schemes for solving the KdV equation, using different approaches to couple the sub-RNNs.

To demonstrate the LRNN-DG methods, we begin by examining a simple linear problem in a one-dimensional spatial domain $\Omega=(x_0, x_{N_s})$. The linear KdV equation is expressed as:
\begin{equation}
\label{linear_kdv}
u_t(t,x) + u_{xxx}(t,x) = 0,\,\, (t,x)\in I\times\Omega,
\end{equation}
where $I$ and $\Omega$ denote the time and space domains, respectively. In constructing the LRNN-DG and LRNN-$C^1$DG methods, we consider the following boundary conditions:
\begin{align}
u(t, x_0) &=  g_0, \label{bc_linear}\\
\frac{\partial u}{\partial x}(t, x_{N_s}) &= g_1, \label{bc_linear2}\\
\frac{\partial^2 u}{\partial x^2}(t, x_{N_s}) &= g_2. \label{bc_linear3}
\end{align}
where $g_i$, $i=0,1,2$, are given functions. The initial condition is specified as:
\begin{align}
\label{ini_linear}
u(t_0, x) = u_0.
\end{align}
where $u_0$ is a known function.

\subsection{LRNN-DG Formulation}

In this subsection, we introduce a space-time DG scheme to integrate the sub-networks effectively.

To derive the LRNN-DG formulation, we start by presenting some identities:
\begin{align}
\int_K \nabla v \cdot \boldsymbol{q}\,\dx &= -\int_K v~(\nabla\cdot \boldsymbol{q}) \,\dx + \int_{\partial K} v\,\boldsymbol{q}\cdot\bn_K\,\ds,\label{ibps}\\
\sum_{K\in \mathcal{T}_h} \int_{\partial K} v \boldsymbol{q} \cdot \boldsymbol{n}_K \, \ds
&= \int_{{\cal E}_h} \llbracket v\rrbracket \cdot \{\boldsymbol{q}\}\, \ds
+ \int_{{\cal E}_h^i} \{v\} \cdot [\boldsymbol{q}]\, \ds,\label{iden_dg}\\
\sum_{I_i\in \mathcal{D}_\tau} (v w)|^{t_i}_{t_{i-1}}
&= \sum^{N_t}_{i=0} \,[v(t_i,\bx)]  \{w(t_i,\bx)\}
+ \sum^{N_t -1}_{i=1} \,\{v(t_i,\bx)\} [w(t_i,\bx)].\label{iden_tdg}
\end{align}

Next, we multiply Equation \eqref{linear_kdv} by a test function $v$ and integrate over the subdomain $\sigma = I_i \times K_j$, where $K_j = (x_{j-1}, x_j)$ represents the one-dimensional spatial subdomain. Integration by parts leads to:
\begin{align*}
&-\int_\sigma u\frac{\partial v}{\partial t}\dt\dx + \int_{K_j} \left({u}v \right)|^{t_i}_{t_{i-1}}\dx - \int_\sigma u \frac{\partial^3 v}{\partial x^3}\dt\dx \\
&+ \int_{I_i} \left.\left( {u} \frac{\partial^2 v}{\partial x^2}\right)\right|^{x_j}_{x_{j-1}} \dt - \int_{I_i} \left. \left(  {\frac{\partial u}{\partial x}} {\frac{\partial v}{\partial x}} \right)\right|^{x_j}_{x_{j-1}} \dt +\int_{I_i}\left. \left( {\frac{\partial^2 u}{\partial x^2}} v \right)\right|^{x_j}_{x_{j-1}} \dt = 0
\end{align*}
We utilize numerical traces such as $\widetilde{u_{h\tau}}$ to approximate $u$ over temporal edges and $\widehat{u_{h\tau}}$, $\widehat{p_{h\tau}}$, and $\widehat{q_{h\tau}}$ to approximate $u$, ${\frac{\partial u}{\partial x}}$, and ${\frac{\partial^2 u}{\partial x^2}}$ over spatial edges. This results in:
\begin{align*}
&-\int_\sigma u \frac{\partial v}{\partial t}\dt\dx + \int_{K_j} \left.\left(\widetilde{u_{h\tau}} v \right)\right|^{t_i}_{t_{i-1}}\dx - \int_\sigma u \frac{\partial^3 v}{\partial x^3}\dt\dx \\
&+ \int_{I_i} \left.\left( \widehat{u_{h\tau}} \frac{\partial^2 v} {\partial x^2}\right)\right|^{x_j}_{x_{j-1}} \dt - \int_{I_i} \left. \left(  {\widehat{p_{h\tau}}} {\frac{\partial v}{\partial x}} \right)\right|^{x_j}_{x_{j-1}} \dt +\int_{I_i}\left. \left( {\widehat{q_{h\tau}}} v \right)\right|^{x_j}_{x_{j-1}} \dt = 0.
\end{align*}

To discretize the continuous problem, we introduce subscripts $h$ and $\tau$, choose a trial function $u_{h\tau} \in V_{h\tau}$ and a test function $v_{h\tau} \in V_{h\tau}$. After adding all terms together and applying the identities \eqref{ibps} and \eqref{iden_tdg}, we arrive at:
\begin{align}
\label{w_kdv_for_flux}
&\int_\Sigma \frac{\partial_\tau u_{h\tau}}{\partial_\tau t} v_{h\tau} \dt\dx + \sum^{N_t}_{i=0}\int_{\mathcal{T}_h} \left[\widetilde{u_{h\tau}}(t_i,x) - u^\tau_h(t_i,x) \right] \left\{v_{h\tau}(t_i,x) \right\} \dx \nonumber\\
&+ \sum^{N_t -1}_{i=1}\int_{\mathcal{T}_h} \left\{\widetilde{u_{h\tau}}(t_i,x) - u^\tau_h(t_i,x) \right\} \left[v_{h\tau}(t_i,x) \right] \dx 
+ \int_{\Sigma} \frac{\partial_h u_{h\tau}}{\partial_h x} \frac{\partial^2_h v_{h\tau}}{\partial_h x^2}\dt\dx\nonumber \\
&+ \sum^{N_s}_{j=0}\int_{\mathcal{D}_h} \left[\widehat{u_{h\tau}}(t, x_j) - u_{h\tau}(t, x_j)\right] \left\{\frac{\partial_h^2 v_{h\tau}}{\partial_h x^2}(t, x_j)\right\} \dt + \sum^{N_e -1}_{j=1}\int_{\mathcal{D}_h} \left\{\widehat{u_{h\tau}}(t, x_j) - u_{h\tau}(t, x_j)\right\} \left[\frac{\partial_h^2 v_{h\tau}}{\partial_h x^2}(t, x_j)\right] \dt\nonumber\\
&- \sum^{N_s}_{j=0}\int_{\mathcal{D}_h}\left\{\widehat{p_{h\tau}}(t, x_j)\right\} \left[{\frac{\partial_h v_{h\tau}}{\partial_h x}(t, x_j)}\right]  \dt - \sum^{N_s -1}_{j=1}\int_{\mathcal{D}_h}\left[\widehat{p_{h\tau}}(t, x_j)\right] \left\{\frac{\partial_h v_{h\tau}}{\partial_h x}(t, x_j)\right\}  \dt\nonumber\\
&+\sum^{N_s}_{j=0}\int_{\mathcal{D}_h} \left\{ {\widehat{q_{h\tau}}(t, x_j)}\right\} \left[ v_{h\tau}(t, x_j) \right] \dt +\sum^{N_s -1}_{j=1}\int_{\mathcal{D}_h} \left[ {\widehat{q_{h\tau}}(t, x_j)}\right] \left\{ v_{h\tau}(t, x_j) \right\} \dt = 0,
\end{align}
where $\frac{\partial_h v_{h\tau}}{\partial_h x}$ represents the broken partial derivative of $v_{h\tau}$ with respect to the mesh $\mathcal{T}_h$ and $\frac{\partial_\tau v_{h\tau}}{\partial_\tau t}$ signifies the broken partial derivative of $v_{h\tau}$ concerning the partition $\mathcal{D}_\tau$.

The numerical fluxes for the specific problem ensure consistency in the LRNN-DG formulation. We select the following fluxes:
\begin{align*}
\widetilde{u_{h\tau}} = \left\{ u_{h\tau}\right\} - \eta_1\left[u_{h\tau}\right] \quad &{\rm on} \quad f \in \mathcal{P}^i_\tau \times \mathcal{T}_h, \\
\widetilde{u_{h\tau}} = u_0 \quad &{\rm on}\quad f\in \{t_0\} \times \mathcal{T}_h,\\
\widetilde{u_{h\tau}} = u_{h\tau} \quad &{\rm on}\quad f\in \{t_{N_t}\} \times \mathcal{T}_h,\\
\widehat{u_{h\tau}} = \left\{ u_{h\tau}\right\}, \widehat{p^\tau_h} = \left\{ \frac{\partial_h u_{h\tau}}{\partial_h x}\right\} - \eta_2\left[u_{h\tau}\right],\widehat{q_{h\tau}} = \left\{ \frac{\partial^2_h u_{h\tau}}{\partial_h x^2}\right\} - \eta_2\left[\frac{\partial_h u_{h\tau}}{\partial_h x}\right] \quad &{\rm on} \quad f \in \mathcal{D}_\tau \times \mathcal{E}^i_h, \\
\widehat{u_{h\tau}} =g_0,\widehat{p_{h\tau}} =\frac{\partial_h u_{h\tau}}{\partial_h x} - \eta_2(g_0 - u_{h\tau}), \widehat{q_{h\tau}} =\frac{\partial^2_h u_{h\tau}}{\partial_h x^2}\quad &{\rm on} \quad f \in \{x_0 \} \times \mathcal{E}^i_h,\\
\widehat{u_{h\tau}} =u_{h\tau},\widehat{p_{h\tau}} = g_1, \widehat{q_{h\tau}} = g_2 - \eta_3(g_1 - \frac{\partial_h u_{h\tau}}{\partial_h x})\quad &{\rm on} \quad f \in \{x_{N_s} \} \times \mathcal{E}^i_h.
\end{align*}
Here, $\eta_1$, $\eta_2$, and $\eta_3$ are penalty parameters set as constants on each edge $f$.

By taking these fluxes into Equation \eqref{w_kdv_for_flux}, we derive the space-time DG formulation of the linear problem as:
\begin{align}
\label{w_kdv_for}
&\int_\Sigma \frac{\partial_\tau u_{h\tau}}{\partial_\tau t} v_{h\tau} \dt\dx - \sum^{N_t -1}_{i=0}\int_{\mathcal{T}_h} \left[ u_{h\tau}(t_i,x) \right] \left\{v_{h\tau}(t_i,x) \right\} \dx - \sum^{N_t -1}_{i=1}\eta_1\int_{\mathcal{T}_h} \left[ u_{h\tau}(t_i,x) \right] \left[v_{h\tau}(t_i,x) \right] \dx \nonumber\\
&+ \int_{\Sigma} \frac{\partial_h u_{h\tau}}{\partial_h x} \frac{\partial^2_h v_{h\tau}}{\partial_h x^2}\dt\dx - \sum^{N_s -1}_{j=0}\int_{\mathcal{D}_h} \left[ u_{h\tau}(t, x_j)\right] \left\{\frac{\partial_h^2 v_{h\tau}}{\partial_h x^2}(t, x_j)\right\} \dt \nonumber \\
&- \sum^{N_s -1}_{j=0}\int_{\mathcal{D}_h}\left\{\frac{\partial_h u_{h\tau}}{\partial_h x}(t, x_j)\right\} \left[{\frac{\partial_h v_{h\tau}}{\partial_h x}(t, x_j)}\right]  \dt + \sum^{N_s -1}_{j=0}\eta_2 \int_{\mathcal{D}_h}\left[ u_{h\tau}(t, x_j)\right]\left[ \frac{\partial_h v_{h\tau}}{\partial_h x}(t, x_j)\right] \dt \nonumber\\
&+ \sum^{N_s -1}_{j=0}\int_{\mathcal{D}_h}\left\{\frac{\partial^2_h u_{h\tau}}{\partial_h x^2}(t, x_j)\right\} \left[{v_{h\tau}}(t, x_j)\right]  \dt \nonumber - \sum^{N_s}_{j=1}\eta_3 \int_{\mathcal{D}_h}\left[ \frac{\partial_h u_{h\tau}}{\partial_h x}(t, x_j)\right]\left[ v_{h\tau}(t, x_j)\right] \dt \nonumber\\
= &- \int_{\mathcal{D}_h}-g_0 \frac{\partial_h^2 v_{h\tau}}{\partial_h x^2}(t, x_0) \dt+ \int_{\mathcal{D}_h} g_1 \frac{\partial_h v_{h\tau}}{\partial_h x}(t, x_{N_s}) \dt + \eta_2  \int_{\mathcal{D}_h}g_0 \frac{\partial_h v_{h\tau}}{\partial_h x}(t, x_0) \dt \nonumber\\
& - \int_{\mathcal{D}_h} g_2  v_{h\tau}(t, x_{N_s}) \dt - \eta_3 \int_{\mathcal{D}_h}g_1 v_{h\tau}(t, x_{N_s}) \dt 
\end{align}

For the nonlinear problem with the term $b^{KdV}(u, \frac{\partial u}{\partial x})$, the equation becomes:
\begin{equation}
\label{nonlinear_kdv}
u_t(t,\bx) + b^{KdV}(u, \frac{\partial u}{\partial x}) + u_{xxx}(t,\bx) = 0,\,\, (t,\bx)\in I\times\Omega,
\end{equation}
subject to the boundary and initial conditions \eqref{bc_linear}--\eqref{ini_linear}. 
The LRNN-DG scheme for solving the nonlinear KdV Equation \eqref{nonlinear_kdv} is: Find $u^\tau_h \in V^\tau_h$ such that
\begin{align}\label{kdvfor}
B^{KdV}({u_{h\tau}}^{(n+1)},v_{h\tau}) + \left(b^{KdV}_L({u_{h\tau}}^{(n+1)}, {u_{h\tau}}^{(n)}),v_{h\tau}\right)\nonumber\\
= l^{KdV}(v_{h\tau})+ \left(b^{KdV}_R({u_{h\tau}}^{(n)}),v_{h\tau}\right) \quad \forall v_{h\tau}\in V^\tau_h\ n=0,1,\cdots.
\end{align}
Here, $B^{KdV}$ and $l^{KdV}$ represent the terms on the left and right sides of Equation \eqref{w_kdv_for}, respectively, while $b^{KdV}_L$ and $b^{KdV}_R$ denote linearizations computed using Newton or Picard methods. The initial value for the iterative algorithm is denoted as ${u_{h\tau}}^{(0)}$.
After a series of nonlinear iterations and least squares computations, numerical solutions can be obtained.

For the nonlinear term $u u_x$, the Newton method yields the following linearization:
$$b^{KdV}_L({u_{h\tau}}^{(n+1)}, {u_{h\tau}}^{(n)}) = \frac{\partial u_{h\tau}^{(n+1)}}{\partial x} {u_{h\tau}}^{(n)} + {u_{h\tau}}^{(n+1)}\frac {\partial u_{h\tau}^{(n)}}{\partial x},$$ 
$$b^{KdV}_R({u_{h\tau}}^{(n)}) = \frac{\partial u_{h\tau}^{(n)}}{\partial x} {u_{h\tau}}^{(n)} .$$
In the case where $b^{KdV}(u, \frac{\partial u}{\partial x}) = u^3 u_x$, the Newton linearization becomes:
$$b^{KdV}_L({u_{h\tau}}^{(n+1)}, {u_{h\tau}}^{(n)}) =3 ({u_{h\tau}}^{(n)})^2 \frac{\partial u_{h\tau}^{(n)}}{\partial x} {u_{h\tau}}^{(n+1)} +   ({u_{h\tau}}^{(n)})^3 \frac {\partial u_{h\tau}^{(n+1)}}{\partial x},$$ 
$$b^{KdV}_R({u_{h\tau}}^{(n)}) = 3 ({u_{h\tau}}^{(n)})^3 \frac{\partial u_{h\tau}^{(n)}}{\partial x}. $$

\begin{remark}
It is important to note that the formulation \eqref{kdvfor} is specifically tailored to solve the KdV Equation \eqref{nonlinear_kdv} with the boundary conditions \eqref{bc_linear}--\eqref{bc_linear3} and the initial condition \eqref{ini_linear}. For different boundary conditions, such as periodic boundaries, distinct flux choices are necessary. Additionally, adjustments to the formulation are required for generalized KdV equations. 
\end{remark}

\subsection{LRNN-$C^1$DG Formulation}

In this subsection, we present an alternative approach that enforces initial conditions, boundary conditions, and continuity conditions at selected collocation points to effectively integrate local neural networks across different subdomains.

Consider the nonlinear problem \eqref{nonlinear_kdv}. By multiplying both sides of the equation by the test function $v_{h\tau}\in V_{h\tau}$ and integrating over a subdomain $\sigma$, we derive, through integration by parts and linearization, the following weak form:
\begin{align}\label{lrnn_c1dg_weak}
{B}^{KdV}_{\sigma}({u_{h\tau}}^{(n+1)},v_{h\tau}) + \left(b^{KdV}_L({u_{h\tau}}^{(n+1)}, {u_{h\tau}}^{(n)}),v_{h\tau}\right)_\sigma= \left(b^{KdV}_R({u_{h\tau}}^{(n)}),v_{h\tau}\right)_\sigma \quad \forall v_{h\tau}\in V_{h\tau},
\end{align}
for $n=0,1,\cdots$,
where
\begin{align*}
{B}^{KdV}_{\sigma}(u_{h\tau},v_{h\tau})=
&- \int_{\sigma} u_{h\tau} \frac{\partial v_{h\tau}}{\partial t} \dt \dx + \int_{K_j} \left({u_{h\tau}}v_{h\tau} \right)|^{t_i}_{t_{i-1}}\dx \\
&- \int_{\sigma}\frac{\partial^2 u_{h\tau}}{\partial x^2} \frac{\partial_h v_{h\tau}}{\partial_h x}\dt\dx + \int_{I_i} \left. \left(  {\frac{\partial^2 u_{h\tau}}{\partial x^2}} {v_{h\tau}} \right)\right|^{x_j}_{x_{j-1}} \dt.
\end{align*}
Where, $\left(\cdot, \cdot\right)_\sigma$ denotes a local inner product within subdomain $\sigma$. Additionally, we introduce supplementary equations to ensure that the numerical solution $u_{h\tau}$ satisfies the initial conditions, boundary conditions and continuity conditions. This leads to the LRNN-$C^1$DG method for the non-linear KdV equation: Find $u_{h\tau} \in V_{h\tau}$ such that $u_{h\tau}$ satisfies Equation \eqref{lrnn_c1dg_weak} in every subdomain $\sigma\in \mathcal{D}_\tau \times \mathcal{T}_h$ and
\begin{align}
u_{h\tau}(t^0_\tau, x^0_\tau) = u_0\quad&\forall (t^0_\tau, x^0_\tau)\in P^0_\tau,\label{c1_t0}\\
[u_{h\tau}(t^i_\tau, x^i_\tau)] = 0 \quad&\forall (t^i_\tau, x^i_\tau)\in P^i_\tau,\label{c1_ti}\\
u_{h\tau}(t^-_h, x^-_h) = g_0\quad&\forall(t^-_h, x^-_h)\in P^-_h,\label{c1_x-}\\
\frac{\partial u_{h\tau}}{\partial x}(t^+_h, x^+_h) = g_1,\frac{\partial^2 u_{h\tau}}{\partial x^2}(t^+_h, x^+_h) = g_2\quad&\forall(t^+_h, x^+_h)\in P^+_h,\label{c1_x+}\\
[u_{h\tau}(t^i_h, x^i_h)] = 0, \left[\frac{\partial u_{h\tau}}{\partial x}(t^i_h, x^i_h)\right], \left[\frac{\partial^2 u_{h\tau}}{\partial x^2}(t^i_h, x^i_h)\right] = 0\quad&\forall (t^i_h, x^i_h)\in P^i_h.\label{c1_xi}
\end{align}
Here, Equation \eqref{c1_t0} imposes the initial condition on the solution $u_{h\tau}$, while Equations \eqref{c1_x-} and \eqref{c1_x+} ensure that $u_{h\tau}$ satisfies the boundary conditions. Equation \eqref{c1_ti} and Equation \eqref{c1_xi} enforce $C^0$ and $C^1$ continuity on temporal edges and spatial edges, respectively. The numerical solution is then obtained using least-squares methods and nonlinear iterations.

\begin{remark}
In \cite{Sun2022lrnndg, Sun2024lrnndgdvwe}, another method called the LRNN-$C^0$DG method combines the LRNN-DG method formulation with certain continuity conditions from the LRNN-$C^1$DG method. In this work, we focus on the LRNN-DG method and the LRNN-$C^1$DG method for solving nonlinear KdV and Burgers equations, as they provide a sufficient illustration of the concepts behind these methods. 
\end{remark}

\section{LRNN-DG Methods for the Burgers Equation}
\label{dgforburgers}

In the previous section, we developed LRNN-DG methods for solving the KdV equation. In this section, we focus on the Burgers equation and present formulations for both the LRNN-DG and the LRNN-$C^1$DG methods.

Consider the Burgers equation with Dirichlet boundary conditions:
\begin{align}
u_t + u(\nabla u\cdot \mathcal{I}) -\epsilon\Delta u = 0 \quad &{\rm in}\quad I\times\Omega,\label{burgers_str}\\
u(t_0,\bx) = u_0,\quad &{\rm on}\quad \{t_0\}\times \Omega,\\
u(t,\bx) = g\quad &{\rm on}\quad I\times \partial \Omega. 
\end{align}
Here, $I=(t_0, T)$ and $\Omega\subset\mathbb{R}^d$ represent the time and space domain, respectively. $\mathcal{I}$ is a d-dimensional vector $[1,1,\cdots,1]^{\rm T}$, and $g$ is a given function. Using the notation introduced in Section \ref{dgforkdv}, we derive the following equation of the linear part through a similar process:
\begin{align}
\label{w_burgers_for_flux}
&\int_\Sigma \frac{\partial_\tau u_{h\tau}}{\partial_\tau t} v_{h\tau} \dt\dx + \sum^{N_t}_{i=0}\int_{\mathcal{T}_h} \left[\widetilde{u_{h\tau}}(t_i,\bx) - u_{h\tau}(t_i,\bx) \right] \left\{v_{h\tau}(t_i,\bx) \right\} \dx \nonumber\\
&+ \sum^{N_t -1}_{i=1}\int_{\mathcal{T}_h} \left\{\widetilde{u_{h\tau}}(t_i,\bx) - u_{h\tau}(t_i,\bx) \right\} \left[v_{h\tau}(t_i,\bx) \right] \dx 
+ \int_{\Sigma}  \nabla_h u_{h\tau} \cdot \nabla_h v_{h\tau}\dt\dx\nonumber \\
&+ \int_{\mathcal{D}_h\times\mathcal{E}_h} \left[\widehat{u_{h\tau}} - u_{h\tau}\right]\cdot \left\{\nabla_h v_{h\tau}\right\} \dt\ds + \int_{\mathcal{D}_h\times\mathcal{E}^i_h} \left\{\widehat{u_{h\tau}} - u_{h\tau}\right\} \cdot\left[\nabla_h v_{h\tau}\right] \dt\ds\nonumber\\
&- \int_{\mathcal{D}_h\times\mathcal{E}_h}\left\{\widehat{\boldsymbol{p}_{h\tau}}\right\} \cdot \left[{v_{h\tau}}\right]  \dt\ds - \int_{\mathcal{D}_h\times\mathcal{E}^i_h}\left[\widehat{\boldsymbol{p}_{h\tau}}\right]\cdot \left\{v_{h\tau}\right\}  \dt\ds = 0.
\end{align}
Here, $\nabla_h$ denotes the broken spatial gradient, $\widehat{\boldsymbol{p}_{h\tau}}\in\boldsymbol{Q}_{h\tau}$ is the approximation of $\nabla_h u_{h\tau}$ on spatial faces, and other terms remain the same as defined in the previous section.

We define the numerical fluxes as follows:
\begin{align*}
\widehat{u_{h\tau}} = \{ u_{h\tau} \} \quad &{\rm on} \ f\in\mathcal{D}_\tau\times\mathcal{E}^i_h,\\
\widehat{u_{h\tau}} =g \quad &{\rm on} \ f\in\mathcal{D}_\tau\times\mathcal{E}^\partial_h,\\
\widetilde{u_{h\tau}} = \{ u_{h\tau} \} - \eta \left[ u_{h\tau} \right] \quad &{\rm on} \ f\in\mathcal{P}^i_\tau\times\mathcal{T}_h,\\
\widetilde{u_{h\tau}} = u_0 \quad &{\rm on} \ f\in\{t_0\}\times\mathcal{T}_h,\\
\widetilde{u_{h\tau}} = u_{h\tau} \quad &{\rm on} \ f\in\{t_{N_t}\}\times\mathcal{T}_h,\\
\widehat{\boldsymbol{p}_{h\tau}} = \{ \nabla_h u_{h\tau} \} - \eta \llbracket u_{h\tau}\rrbracket \quad &{\rm on} \ f\in \mathcal{D}_\tau\times\mathcal{E}^i_h,\\
\widehat{\boldsymbol{p}_{h\tau}} = \nabla_h u_{h\tau} - \eta (u-g)\boldsymbol{n} \quad &{\rm on} \ f\in \mathcal{D}_\tau\times\mathcal{E}^\partial_h,
\end{align*}
where $\eta$ is a penalty parameter, constant on each face $f$ that ensures the numerical solution $u_{h\tau}$ satisfies initial condition, boundary condition, and continuity.

Thus, we obtain the LRNN-DG method for the Burgers equation: Find $u_{h\tau}\in V_{h\tau}$, such that
\begin{align}\label{burgers_dgfor}
B^{Bur}({u_{h\tau}}^{(n+1)},v_{h\tau}) + \left(b^{Bur}_L({u_{h\tau}}^{(n+1)}, {u_{h\tau}}^{(n)}),v_{h\tau}\right) = l^{Bur}(v_{h\tau})+ \left(b^{Bur}_R({u_{h\tau}}^{(n)}),v_{h\tau}\right) \quad \forall v_{h\tau}\in V_{h\tau}
\end{align}
for $n=0,1,\cdots$,
where $\left(b^{Bur}_L({u_{h\tau}}^{(n+1)}, {u_{h\tau}}^{(n)}),v_{h\tau}\right)$ and $\left(b^{Bur}_R({u_{h\tau}}^{(n)}),v_{h\tau}\right)$ represent the linearizations of the nonlinear term $u(\nabla u\cdot \mathcal{I})$ and
\begin{align}
B^{Bur}(u_{h\tau},v_{h\tau}) &= 
\int_{\Sigma} \frac{\partial_\tau u_{h\tau}}{\partial_\tau t}v_{h\tau} \dt \dx + \int_{\Sigma} \nabla_h u_{h\tau} \cdot\nabla_h v_{h\tau} \dt \dx\nonumber\\ 
&\quad-\sum^{N_t -1}_{i=0}\int_{\mathcal{T}_h} \left[ u_{h\tau}(t_i,\bx)\right] \cdot \{ v_{h\tau}(t_i,\bx) \} \dx-\sum^{N_t -1}_{i=1} \int_{\mathcal{T}_h} \eta \left[ u_{h\tau}(t_i,\bx) \right] \cdot \left[ v_{h\tau}(t_i,\bx)  \right] \dx \nonumber\\
&\quad- \int_{\mathcal{D}_\tau \times\mathcal{E}^{}_h} \left(\llbracket u_{h\tau}\rrbracket\cdot \{ \nabla_{h}v_{h\tau} \} +\llbracket v_{h\tau}\rrbracket\cdot \{ \nabla_{h}u_{h\tau} \} -\eta \llbracket u_{h\tau}\rrbracket\cdot \llbracket v_{h\tau}\rrbracket\right)\dt \ds,\\
l^{Bur}(v_{h\tau})
&= \int_{\Sigma} f v_{h\tau} \dt \ds
- \int_{ \mathcal{D}_\tau\times{\cal E}_h^{\partial}} (g \boldsymbol{n}\cdot \nabla_{h} v_{h\tau} - \eta g v_{h\tau}) \dt\ds 
+ \int_{\mathcal{T}_h} u_0(\bx) v_{h\tau}(t_0,\bx) \dx.
\end{align}
Here, $\eta = {\eta_f}{(h_f)}^{-1}$, and $\eta_f$ may vary depending on the choice of face $f$. 

Finally, we obtain the numerical solution through nonlinear iterations and solving least-squares problems.

Next, we develop the LRNN-$C^1$DG method for the Burgers equation. Similarly, we multiply both sides of Equation \eqref{burgers_str} by a test function $v_{h\tau} \in V_{h\tau}$ and integrate over a local subdomain $\sigma = I_i\times K$:
\begin{align}\label{lrnn_c1dg_weak_burgers}
{B}^{Bur}_{\sigma}({u_{h\tau}}^{(n+1)},v_{h\tau}) + \left(b^{Bur}_L({u_{h\tau}}^{(n+1)}, {u_{h\tau}}^{(n)}),v_{h\tau}\right)_\sigma =\left(b^{Bur}_R({u_{h\tau}}^{(n)}),v_{h\tau}\right)_\sigma \quad \forall v_{h\tau}\in V_{h\tau}
\end{align}
for $n=0,1,\cdots$,
where
\begin{align*}
{B}^{Bur}_{\sigma}(u_{h\tau},v_{h\tau})=
&- \int_{\sigma} u_{h\tau} \frac{\partial v_{h\tau}}{\partial t} \dt \dx + \int_{K} \left({u_{h\tau}}v_{h\tau} \right)|^{t_i}_{t_{i-1}}\dx \\
&+\int_{\sigma} \nabla u_{h\tau} \cdot \nabla v_{h\tau} \dt\dx - \int_{I_i\times \partial K} \nabla u_{h\tau} \cdot \boldsymbol{n} v_{h\tau} \ds\dt.
\end{align*}
By introducing collocation points and imposing appropriate conditions at these points, we add the following system of equations:
\begin{align}
u_{h\tau}(t^0_\tau, x^0_\tau) = u_0\quad&\forall (t^0_\tau, x^0_\tau)\in P^0_\tau,\label{bur_c1_t0}\\
[u_{h\tau}(t^i_\tau, x^i_\tau)] = 0 \quad&\forall (t^i_\tau, x^i_\tau)\in P^i_\tau,\label{bur_c1_ti}\\
u_{h\tau}(t^\partial_h, x^\partial_h) = g \quad&\forall(t^\partial_h, x^\partial_h)\in P^\partial_h,\label{bur_c1_xb}\\
[u_{h\tau}(t^i_h, x^i_h)] = 0, \left[\nabla u_{h\tau} (t^i_h, x^i_h)\right] = 0\quad&\forall (t^i_h, x^i_h)\in P^i_h.\label{bur_c1_xi}
\end{align}
Thus, we have the LRNN-$C^1$DG scheme for the Burgers equation: Find $u_{h\tau} \in V_{h\tau}$ such that $u_{h\tau}$ satisfies Equation \eqref{lrnn_c1dg_weak_burgers} in each subdomain $\sigma\in \mathcal{D}_\tau\times\mathcal{T}_h$ and meets the strong conditions \eqref{bur_c1_t0}-\eqref{bur_c1_xi}.

\section{Adaptive Domain Decomposition}
\label{sec:adap}

To enhance the accuracy and effectiveness of the proposed method, this sesction introduces adaptive domain decomposition (\cite{Babuska1978adap,Wang2019adap, Wang2015adap}), in which the domain decomposition is guided by error estimators that reflect the error distribution of the numerical solution.

Consider a general PDE of the form
$$
\mathcal{L} u = f  \quad {\rm in}\quad \Sigma,
$$
where $\mathcal{L}$ is a partial diferential operator.
Initially, we generate a uniform rectangular mesh $\mathcal{M}^{(0)}_{\tau h}$ and calculate the interior residuals as follows:
$$\mathcal{R}_{\sigma_i} = \mathcal{L} u - f \quad \forall \sigma_i \in \mathcal{M}^{(0)}_{\tau h}.$$ 
Next, we compute the local error estimators $\{\mathcal{R}^{h}_i = \left(h^2_{\sigma_i} \Vert \mathcal{R}_{\sigma_i}\Vert_{\sigma_i} \right)^{1/2}, i= 1,2,\cdots, N_e\}$ in the $L^2$ norm for each subdomain $\sigma_i\in\mathcal{M}^{(0)}_{\tau h}$. The sum of residuals, $\mathcal{R}^{tol}=\sum_{i=1}^{N_e} (\mathcal{R}_{i}^h)^2$ is then calculated. The local residuals are rearranged from largest to smallest $\{\mathcal{R}_{i}^h, i= 1,2,\cdots, N_e \}$. Additionally, a parameter $\beta$ is introduced as a positive constant. We determine the smallest $N_r$ such that $\sum_{i=1}^{N_r} (\mathcal{R}_{i}^h)^2 \geq \beta\mathcal{R}^{tol}$. Finally, we refine the first $N_r$ subdomains with the lagest local error estimator, leading to a new decomposition $\mathcal{M}^{(1)}_{\tau h}$. This process is repeated until the desired level of accuracy is achieved or the maximum number of iterations is reached.

Alternatively, if we know the characteristic direction of the solution in advance through methods such as traveling wave analysis (\cite{Mohyud2009twkdv1, Kudryashov2009twkdv2, Seadawy2017twkdv3}) or other sources of information, the use of LRNN-DG methods can be more effective for solving KdV-type equations. For example, if the characteristic direction is $x=kt$, we can derive a priori information about the solution's shape from known initial conditions. We then partition the domain $\Sigma$ into different subdomains along the characteristic direction based on this prior information, and introduce distinct local networks to approximate the solution in each subdomain.

Additionally, we can incorporate wavelet basis functions (\cite{Mallat1999wavelet}), aligning them with the inherent characteristic directional information of the problem, to construct a corresponding set of activation functions. This tailored approach enhances the network’s approximation capability for the specific problem, thereby achieving highly accurate numerical results efficiently. Further details on the implementation of these ideas, involving the characteristic direction of the solution, will be demonstrated through numerical examples.

\section{Numerical Examples}
\label{sec:ex}

In this section, we present experiments to demonstrate the effectiveness of the proposed methods in solving the KdV and Burgers equations. 

We begin by introducing some notation. The term ${\rm DoF}_\sigma$ denotes the number of degrees of freedom (DoF) in each subdomain $\sigma$. The global $L^2$ error is defined as:
$$E^{L^2} = \left(\int_\Sigma (u_{h\tau}-u^*)^2 \dx\dt\right)^{\frac{1}{2}},$$ 
and the global $H^1$ error is defined as:
$$E^{H^1} = \left(\int_\Sigma \left(\frac{\partial u_{h\tau}}{\partial t} -\frac{\partial u^*}{\partial t}\right)^2 + |\nabla (u_{h\tau} - u^*)|^2 \dx\dt\right)^{\frac{1}{2}},$$ 
where $u_{h\tau}$ is the numerical solution and $u^*$ represents the exact solution.

We define the difference between the solutions $u_{h\tau}^{(n-1)}$ and $u_{h\tau}^{(n)}$ of two consecutive steps in the nonlinear iteration as:
$$D(u_{h\tau}^{(n)}, u_{h\tau}^{(n-1)}) = \left(\int_\Sigma \left(u_{h\tau}^{(n)} - u_{h\tau}^{(n-1)} \right)^2\dx\dt\right)^{1/2}.$$ 
The stopping criterion of the nonlinear iteration is $D(u_{h\tau}^{(n)}, u_{h\tau}^{(n-1)}) < \epsilon_0$. Additionally, a maximum number $N_{ni}$ of iterations is set to avoid an infinite loop.

We utilize the Pytorch library in Python to construct local neural networks in different subdomains. The Tanh function serves as the activation function, and the parameters $\theta_j$ of hidden layers are randomly generated from a uniform distribution $U(-r,r)$ and remain fixed throughout the training process, where $r$ is a positive constant. The influence of the parameter $r$ is discussed in \cite{Dong2021locELMW}, and strategies for initialization can be found in relevant literature (\cite{Zhang2024Transferable, Dang2024adaptive}). Gaussian quadrature is applied to evaluate all integrals in the experiments. Finally, the least-squares method, implemented using the Scipy package, is used to optimize the output layer parameters. In the numerical examples, we set a fixed random seed to ensure reproducibility.

\begin{example}[Generalized KdV Equation]
\label{ex_gKdV}

In this experiment, we evaluate the performance of the proposed methods by solving the generalized KdV equation (\cite{Yan2002dgkdv}) with a small coefficient for the third derivative term, which features a soliton solution. The equation is given by:
\begin{equation}
u_t + u_x + u^3 u_x +\epsilon u_{xxx} = 0\quad (t,x)\in I\times\Omega,
\end{equation}
where $\Omega = [-2,3]$, $I=[0,2]$ and $\epsilon = 0.2058\times 10^{-4}$. The exact solution is:
\begin{equation}
u(t, x) = A {\rm sech}^{\frac{2}{3}} (K(x-x_0)- \omega_t),
\end{equation}
where $A = 0.2275$, $K = 3(\frac{A^3}{40\epsilon})^{\frac{1}{2}}$, $\omega = K(1+\frac{A^3}{10})$ and $x_0 = 0.5$. The boundary and initial conditions are specified as:
\begin{align}
u(t, -2) = g_0,\quad u_x(t, 3) = g_1,\quad u_{xx}(t,3) = g_2,\quad u(0, x) = u_0,
\end{align}
where $g_0$, $g_1$, $g_2$ and $u_0$ are computed based on the exact solution.
\end{example}

We construct local randomized neural networks on a uniform mesh and compute the global $L^2$ error $E^{L^2} (u)$ and global $H^1$ error $E^{H^1} (u)$ for both the LRNN-DG and LRNN-$C^1$DG methods. For nonlinear iterations, we set $\epsilon = 10^{-6}$ and $N_{ni} = 10$.

\begin{table}[h]
\centering
\begin{tabular}{|c|ll|ll|ll|}
\hline
$\tau,\ h$                   & \multicolumn{2}{c|}{1, 1}                                                            & \multicolumn{2}{c|}{1/2, 1/2}                                 & \multicolumn{2}{c|}{1/4, 1/4}                                 \\ \hline
\diagbox[width=5em,trim=l]{${\rm DoF}_\sigma$}{Norm} & \multicolumn{1}{c|}{$E^{L^{2}}$}                         & \multicolumn{1}{c|}{$E^{H^{1}}$} & \multicolumn{1}{c|}{$E^{L^{2}}$}  & \multicolumn{1}{c|}{$E^{H^{1}}$} & \multicolumn{1}{c|}{$E^{L^{2}}$}  & \multicolumn{1}{c|}{$E^{H^{1}}$} \\ \hline
80                         & \multicolumn{1}{l|}{3.67E-02}                        & 6.32E-01                     & \multicolumn{1}{l|}{3.91E-03} & 1.49E-01 & \multicolumn{1}{l|}{7.38E-05} & 6.84E-03                     \\ \hline
160                        & \multicolumn{1}{l|}{1.03E-02}                        & 2.92E-01                     & \multicolumn{1}{l|}{4.24E-04} & 2.43E-02                     & \multicolumn{1}{l|}{1.00E-05} & 1.23E-03                     \\ \hline
320                        & \multicolumn{1}{l|}{5.47E-03}                        & 1.93E-01                     & \multicolumn{1}{l|}{1.80E-04} & 1.13E-02                     & \multicolumn{1}{l|}{2.14E-06} & 2.76E-04                     \\ \hline
\end{tabular}
\caption{Global errors of the space-time LRNN-DG method in Example \ref{ex_gKdV}}
\label{tablegKdVDG}
\end{table}

\begin{table}[h]
\centering
\begin{tabular}{|c|ll|ll|ll|}
\hline
$\tau, h$                   & \multicolumn{2}{c|}{1, 1}                                                            & \multicolumn{2}{c|}{1/2, 1/2}                                 & \multicolumn{2}{c|}{1/4, 1/4}                                 \\ \hline
\diagbox[width=5em,trim=l]{${\rm DoF}_\sigma$}{Norm} & \multicolumn{1}{c|}{$E^{L^{2}}$}                         & \multicolumn{1}{c|}{$E^{H^{1}}$} & \multicolumn{1}{c|}{$E^{L^{2}}$}  & \multicolumn{1}{c|}{$E^{H^{1}}$} & \multicolumn{1}{c|}{$E^{L^{2}}$}  & \multicolumn{1}{c|}{$E^{H^{1}}$} \\ \hline
80                         & \multicolumn{1}{l|}{ 6.55E-02}                       & 1.24E+00                     & \multicolumn{1}{l|}{8.34E-03} & 2.83E-01                     & \multicolumn{1}{l|}{4.53E-04} & 2.80E-02                     \\ \hline
160                        & \multicolumn{1}{l|}{4.65E-02}                        & 1.32E+00                     & \multicolumn{1}{l|}{4.76E-04} & 2.30E-02                     & \multicolumn{1}{l|}{1.63E-05} & 1.16E-03                     \\ \hline
320                        & \multicolumn{1}{l|}{3.14E-02}                        & 6.92E-01                     & \multicolumn{1}{l|}{1.78E-04} & 8.77E-03                     & \multicolumn{1}{l|}{1.11E-06} & 9.22E-05                     \\ \hline
\end{tabular}
\caption{Global errors of the space-time LRNN-$C^1$DG method in Example \ref{ex_gKdV}}
\label{tablegKdVc1DG}
\end{table}

Table \ref{tablegKdVDG} shows the performance of the LRNN-DG method, presenting errors for various spatial mesh sizes $h$, temporal interval lengths $\tau$, and degrees of freedom in each subdomain. Here, we set $r = 1.76$ and the interior penalty $\eta = 220/h$ or $\eta = 220/\tau$, using 15 Gaussian quadrature points in each direction. The errors decrease as the degrees of freedom increase, and as $h$ and $\tau$ decrease.

Similarly, Table \ref{tablegKdVc1DG} details the errors for the LRNN-$C^1$DG method. Using $r = 1.9$, 17 Gaussian integration points in each direction, and 13 collocation points along each edge, the LRNN-$C^1$DG method demonstrates comparable accuracy to the LRNN-DG method.

\begin{figure}[H] 
  \centering  
  \begin{subfigure}{0.4\textwidth}
     \centering      
      \includegraphics[width=1\textwidth]{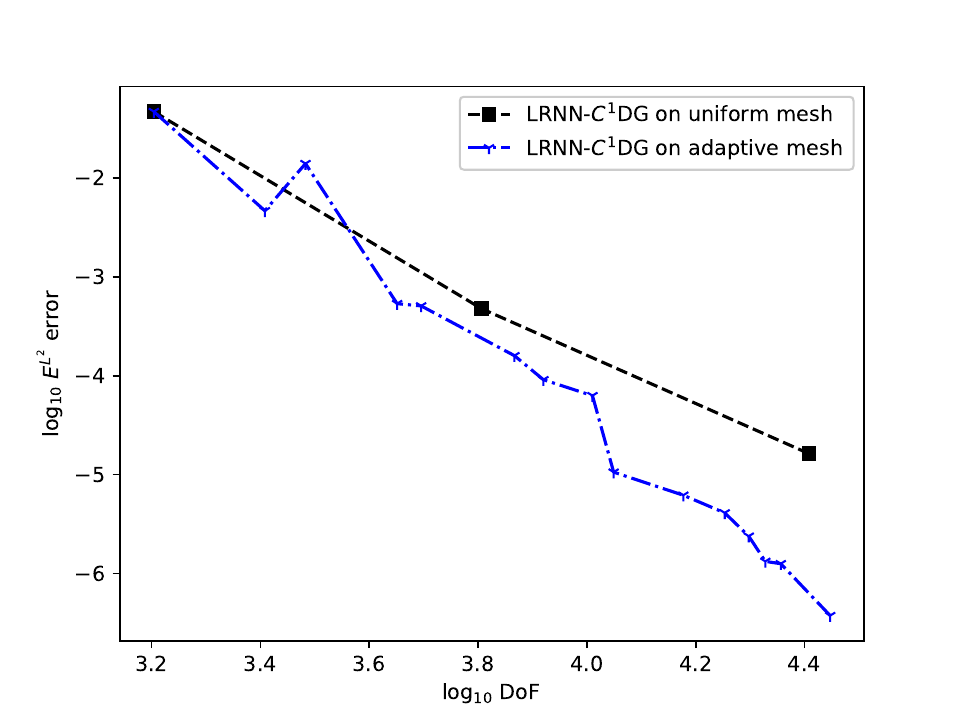}  
      \caption{$E^{L^2}$ errors with respect to DoF}
  \end{subfigure}
  \quad
  \begin{subfigure}{0.4\textwidth}
      \centering      
      \includegraphics[width=1\textwidth]{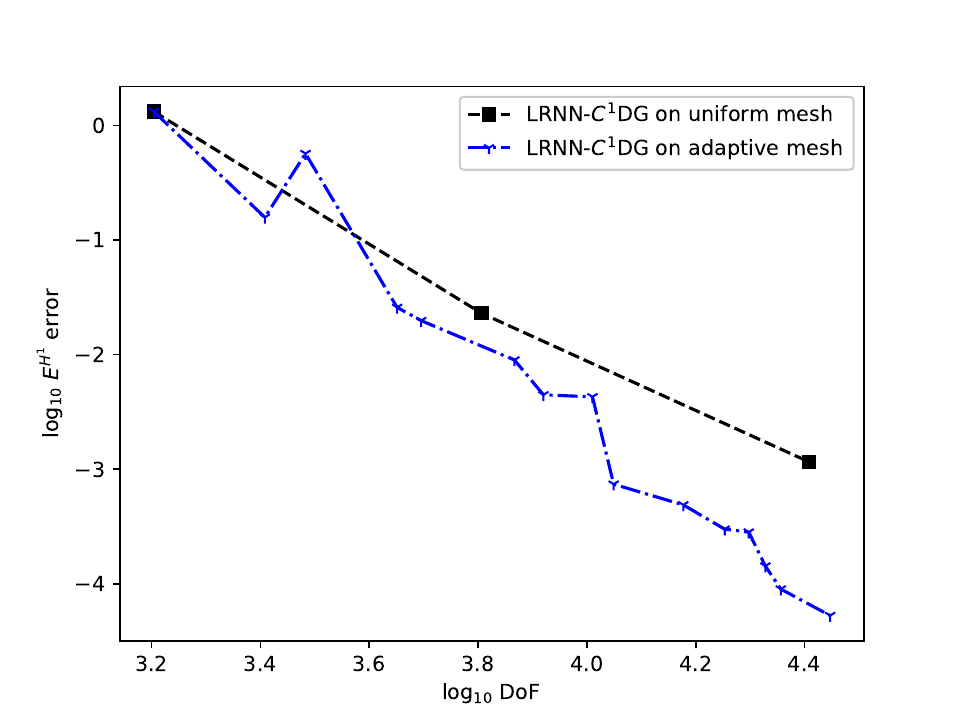}  
      \caption{$E^{H^1}$ errors with respect to DoF}
  \end{subfigure}
  \caption{Errors of LRNN-$C^1$DG methods on uniform and adaptive meshes in Example \ref{ex_gKdV}.}  
  \label{figure_kdv_compare}
\end{figure}

\begin{figure}[H] 
  %\centering  
  \begin{subfigure}{0.24\textwidth}
      %\centering      
      \includegraphics[width=1.1\textwidth]{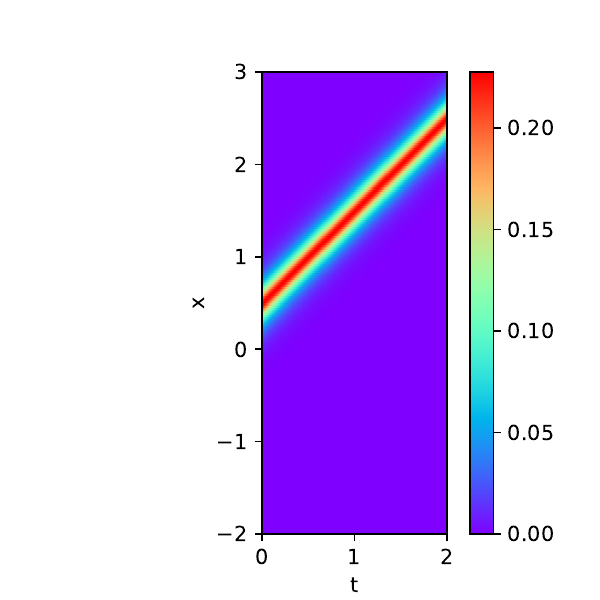}  
      \caption{Exact solution}
      \label{fig:a}
  \end{subfigure}
  \begin{subfigure}{0.24\textwidth}
      %\centering      
      \includegraphics[width=1.1\textwidth]{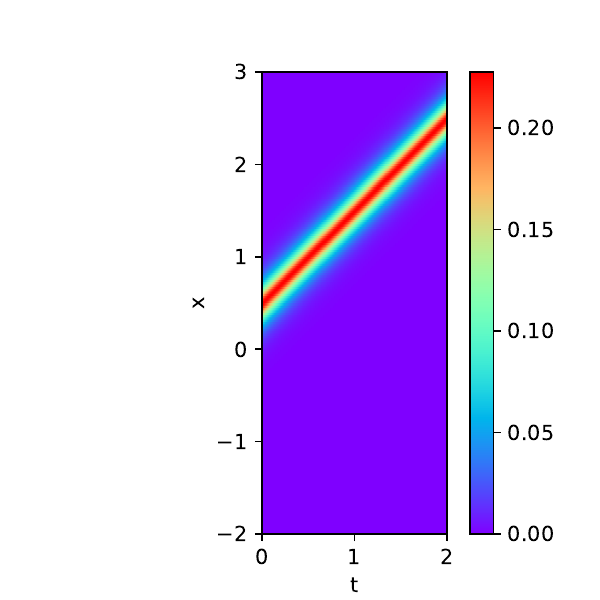}  
      \caption{Numerical solution}
      \label{fig:b}
  \end{subfigure}
  \begin{subfigure}{0.24\textwidth}
      %\centering    
      \includegraphics[width=1.1\textwidth]{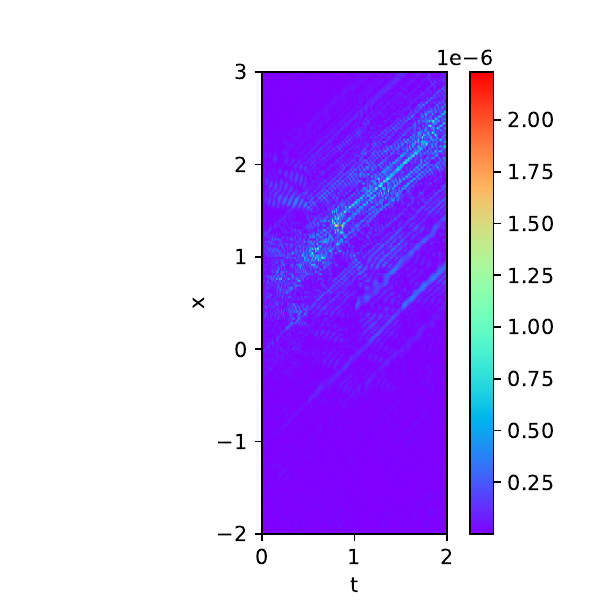}  
		 \caption{$|u_{h\tau} - u^*|$}
      \label{fig:c}
    \end{subfigure}  
   \begin{subfigure}{0.24\textwidth}
      %\centering    
      \includegraphics[width=1.4\textwidth]{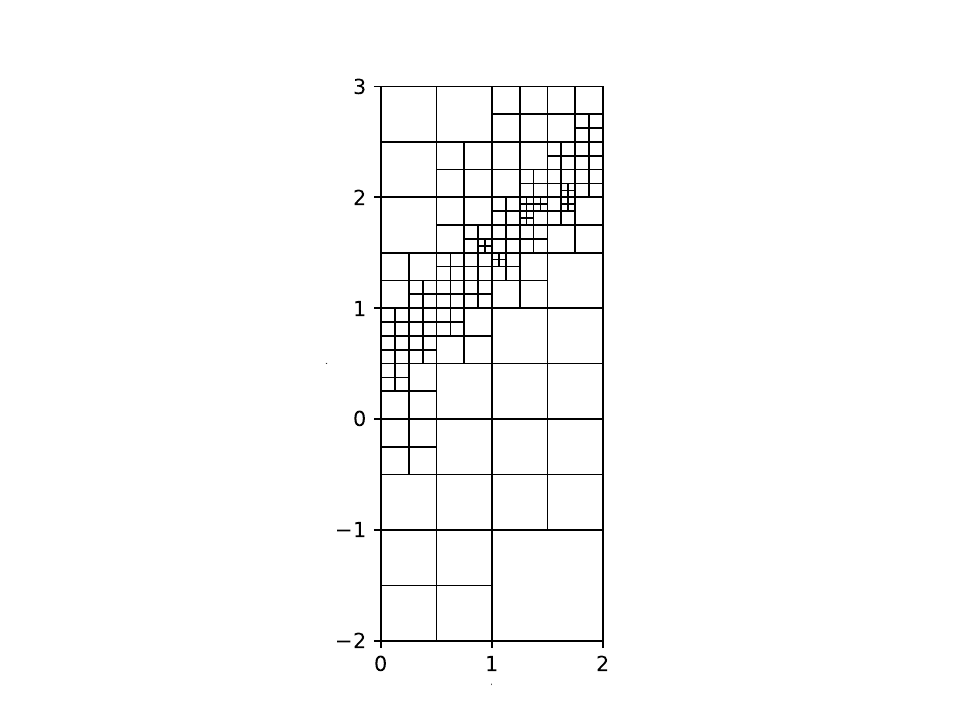}  
		 \caption{Adaptive mesh}
      \label{fig:d}
    \end{subfigure}  
  \caption{The performances of the adaptive LRNN-$C^1$DG method in Example \ref{ex_gKdV}.}  
  \label{figure_adapc1dg}
\end{figure}

Next, we explore the construction of local neural networks on an adaptive mesh, starting from an initial uniform mesh. By setting $\beta = 0.7$, initial $h = 1$, $\tau = 1$, and ${\rm DoF}_\sigma = 160$, we compute errors for the LRNN-$C^1$DG method on the adaptive mesh, as shown in Tables \ref{tablegKdVadap}. To further illustrate the effectiveness of the adaptive meshing technique, Figure \ref{figure_kdv_compare} provides a visual comparison. Additionally, Figures \ref{figure_adapc1dg} display visual comparisons of the exact solution $u^*$, numerical solution $u_{h\tau}$, absolute error, and the adaptive mesh for the proposed methods, showing the success of the adaptive strategy. Figure \ref{figure_adapc1dg}(d) includes 148 subdomains.

%Next, we explore the construction of local neural networks on an adaptive mesh, starting from an initial uniform mesh. By setting $\beta = 0.5$, initial $h=1$, $\tau=1$, and ${\rm DoF}_\sigma = 160$, we record the errors of  LRNN-$C^1$DG methods on the adaptive mesh, as presented in Tables \ref{tablegKdVadap}. To further illustrate the effectiveness of the adaptive meshing technique, we provide Figure \ref{figure_kdv_compare}. Additionally, Figures \ref{figure_adapc1dg} illustrate visual comparisons of the exact solution $u^*$, numerical solution $u_{h\tau}$, absolute error, and the adaptive mesh for the proposed methods, highlighting the success of the adaptive strategy. Figure \ref{figure_adapc1dg} (d) includes 175 subdomains.

\begin{table}[H]
\centering
\begin{tabular}{cccccccccc}
\toprule
{$N_e$} & 10       & 16       & 25       & 37       & 61       & 76       & 88       & 133 & 148       \\ \midrule
$E^{L^2}$     & 4.65E-02  & 6.30E-03. & 8.52E-04 & 5.24E-04 & 2.55E-04 & 5.75E-05 & 3.75E-06 & 2.57E-06 &  9.51E-07 \\
$E^{H^1}$     & 1.32E+00 & 2.20E-01 & 5.33E-02 & 2.73E-02 & 1.14E-02 & 4.64E-03 & 3.79E-04 & 2.28E-04 &  8.27E-05 \\ \bottomrule
\end{tabular}
\caption{Global errors of LRNN-$C^1$DG method on adaptive meshes in Example \ref{ex_gKdV}}
\label{tablegKdVadap}
\end{table}

\vspace{-2mm}

\begin{figure}[H] 
  \centering  
   
  \begin{subfigure}{0.24\textwidth}
      \centering      
      \includegraphics[width=1.1\textwidth]{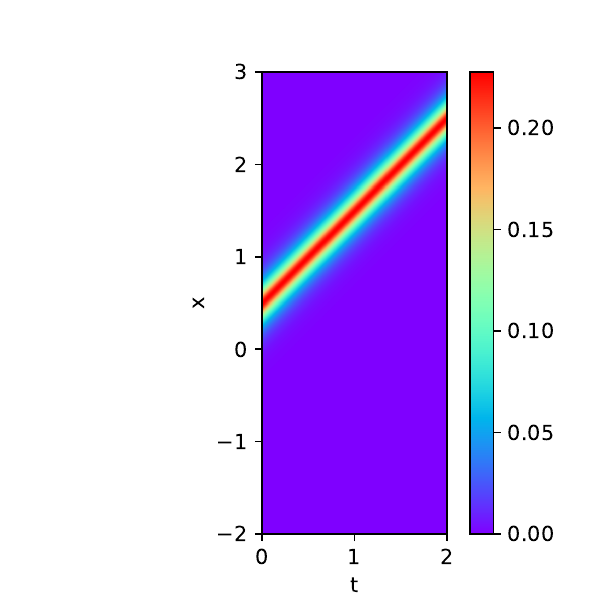}  
      \caption{Exact solution}
      \label{fig:a}
  \end{subfigure}
  \begin{subfigure}{0.24\textwidth}
      \centering      
      \includegraphics[width=1.1\textwidth]{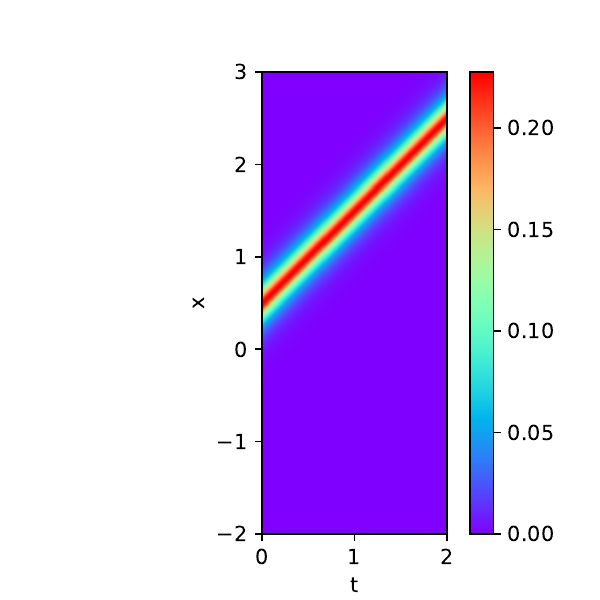}  
      \caption{Numerical solution}
      \label{fig:b}
  \end{subfigure}
  \begin{subfigure}{0.24\textwidth}
      \centering    
      \includegraphics[width=1.1\textwidth]{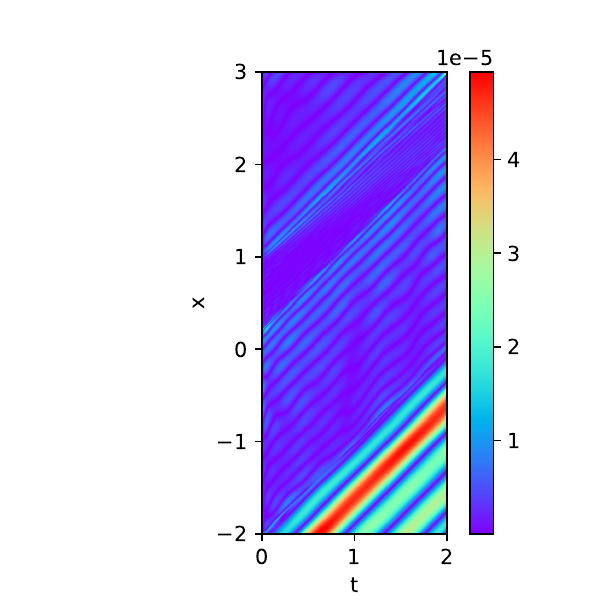}  
		 \caption{$|u_{h\tau} - u^*|$}
      \label{fig:c}
    \end{subfigure}  
   \begin{subfigure}{0.24\textwidth}
      \centering    
      \includegraphics[width=1.1\textwidth]{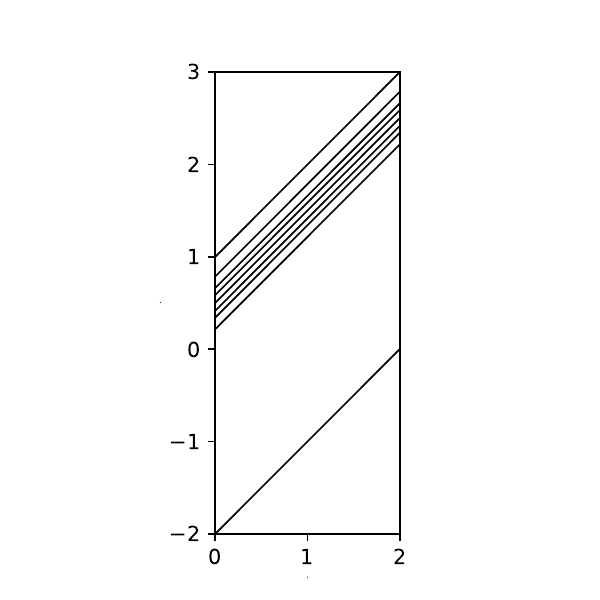}  
		 \caption{\label{spe_mesh}Characteristic mesh}
      %\label{fig:d}
    \end{subfigure}  
  \caption{Performances of the LRNN-DG method on the characteristic mesh in Example \ref{ex_gKdV}.}  
  \label{figure_spedg}
\end{figure}

\begin{figure}[H] 
  \centering  
   
  \begin{subfigure}{0.24\textwidth}
      \centering      
      \includegraphics[width=1.1\textwidth]{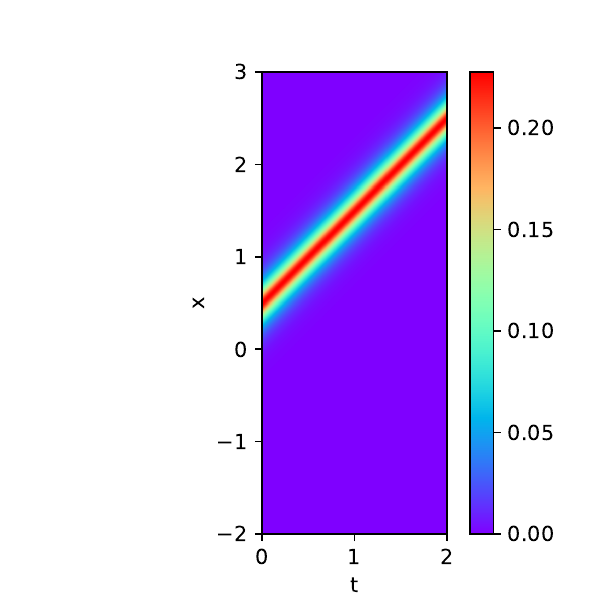}  
      \caption{Exact solution}
      \label{fig:a}
  \end{subfigure}
  \begin{subfigure}{0.24\textwidth}
      \centering      
      \includegraphics[width=1.1\textwidth]{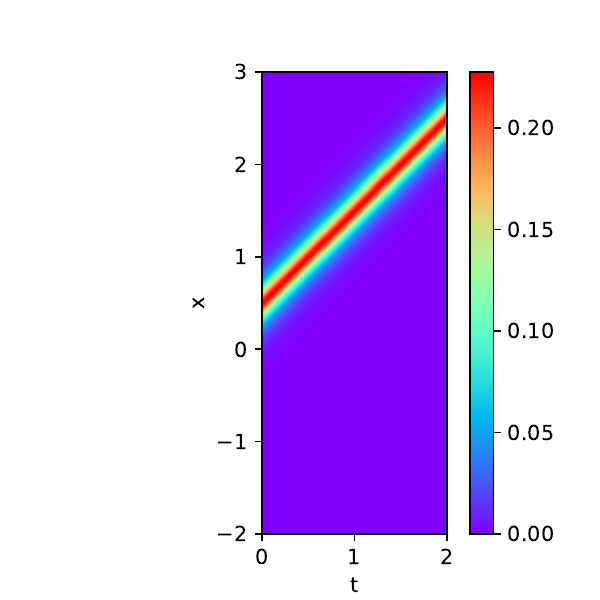}  
      \caption{Numerical solution}
      \label{fig:b}
  \end{subfigure}
  \begin{subfigure}{0.24\textwidth}
      \centering    
      \includegraphics[width=1.1\textwidth]{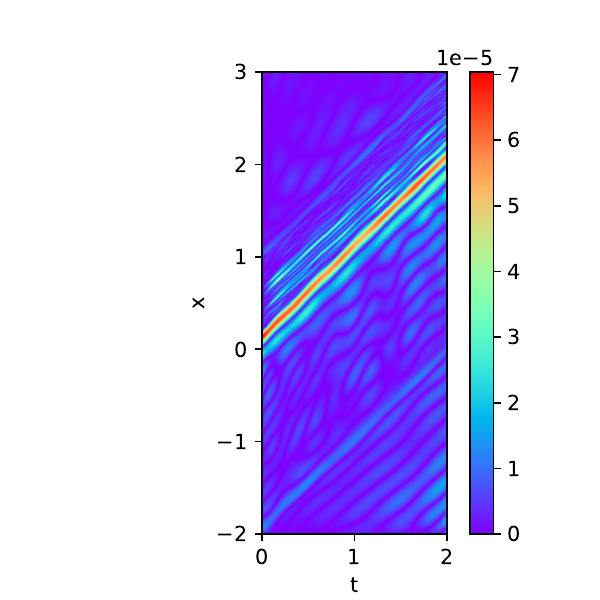}  
		 \caption{$|u_{h\tau} - u^*|$}
      \label{fig:c}
    \end{subfigure}  
   \begin{subfigure}{0.24\textwidth}
      \centering    
      \includegraphics[width=1.1\textwidth]{special_mesh.eps}  
		 \caption{Characteristic mesh}
      \label{fig:d}
    \end{subfigure}  
  \caption{Performances of the LRNN-$C^1$DG method on the characteristic mesh in Example \ref{ex_gKdV}.}  
  \label{figure_spec1dg}
\end{figure}

To enhance both accuracy and efficiency, we also develop a characteristic mesh based on the characteristic direction of the equation, which in this case is given by $t = kx$, where $k = 1.0011774546$ (\cite{Liu2019k}). We identify the set of critical points $P^c_0=\{x: u_0(x)=0\ {\rm or}\ u_{0_{xx}}(x)=0\ {\rm or}\  u_{0_{xxxx}}(x)=0\}$ with $|P^c_0| = 7$. Two additional points are added for simplying numerical integration, resulting in the final mesh depicted in Figure \ref{spe_mesh}. This tailored mesh, with a small number of subdomains, ensures high accuracy, as demonstrated in Table \ref{tablegKdVspe}. Using 17 Gaussian quadrature points in each direction for integration over a parallelogram domain and dividing triangles into 27 smaller sub-triangles with 27 integration points, we achieve high numerical accuracy. For the LRNN-DG method, parameters are set to $r = 1.9$ and $\eta = 200/S_{h\tau}$, where $S_{h\tau}$ represents the edge size. For the LRNN-$C^1$DG method, parameters include $r = 1.9$ and 13 collocation points along each edge. Figures \ref{figure_spedg} and \ref{figure_spec1dg} visually display the exact solutions, numerical solutions, absolute errors, and characteristic mesh decompositions with ${\rm DoF} = 320$. Compared to the uniform and adaptive meshes, the proposed methods achieve significantly higher accuracy with the specialized mesh, using only 10 subdomains.

\begin{table}[h]
\centering
\begin{tabular}{|c|ll|ll|}
\hline
Method                     & \multicolumn{2}{c|}{LRNN-DG}                                                        & \multicolumn{2}{c|}{LRNN-$C^1$DG}                            \\ \hline
\diagbox[width=5em,trim=l]{${\rm DoF}_\sigma$}{Norm} & \multicolumn{1}{c|}{$E^{L^{2}}$}                         & \multicolumn{1}{c|}{$E^{H^{1}}$} & \multicolumn{1}{c|}{$E^{L^{2}}$}  & \multicolumn{1}{c|}{$E^{H^{1}}$} \\ \hline
40                         & \multicolumn{1}{l|}{9.45E-04}                        & 6.17E-02                     & \multicolumn{1}{l|}{1.81E-02} & 2.15E-01                     \\ \hline
80                         & \multicolumn{1}{l|}{8.87E-04} 				  & 3.56E-02                     & \multicolumn{1}{l|}{1.89E-02} & 3.26E-01                     \\ \hline
160                        & \multicolumn{1}{l|}{8.27E-05}                        & 2.19E-03                     & \multicolumn{1}{l|}{8.44E-05} & 1.81E-03                     \\ \hline
320                        & \multicolumn{1}{l|}{4.69E-05}                        & 1.29E-03                     & \multicolumn{1}{l|}{3.83E-05} & 1.31E-03                     \\ \hline
\end{tabular}
\caption{Global errors of LRNN-DG methods on the characteristic mesh in Example \ref{ex_gKdV}}
\label{tablegKdVspe}
\end{table}

\begin{table}[h]
\centering
\begin{tabular}{|c|ll|ll|}
\hline
Method                   & \multicolumn{2}{c|}{LRNN-DG}                                    & \multicolumn{2}{c|}{LRNN-$C^1$DG}                               \\ \hline
\diagbox[width=5em,trim=l]{$DoF_\sigma$}{Norm} & \multicolumn{1}{c|}{$E^{L^2}$} & \multicolumn{1}{c|}{$E^{H^1}$} & \multicolumn{1}{c|}{$E^{L^2}$} & \multicolumn{1}{c|}{$E^{H^1}$} \\ \hline
20                       & \multicolumn{1}{l|}{3.16E-05}  & 2.44E-03                       & \multicolumn{1}{l|}{7.82E-04}  & 2.53E-02                       \\ \hline
40                       & \multicolumn{1}{l|}{1.43E-05}  & 4.16E-04                       & \multicolumn{1}{l|}{3.95E-06}  & 2.55E-04                       \\ \hline
80                       & \multicolumn{1}{l|}{3.90E-06}  & 3.44E-04                       & \multicolumn{1}{l|}{1.60E-06}  & 1.29E-04                       \\ \hline
\end{tabular}
\caption{Global errors of LRNN-DG methods with the wavelet activation function in Example \ref{ex_gKdV}}
\label{tablegKdVwavelet}
\end{table}

\begin{figure}[H] 
  \centering  
   
  \begin{subfigure}{0.24\textwidth}
      \centering      
      \includegraphics[width=1.1\textwidth]{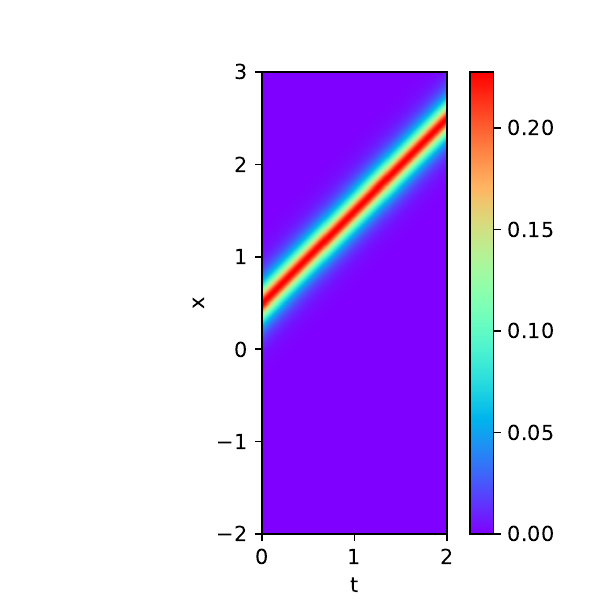}  
      \caption{Exact solution}
      \label{fig:a}
  \end{subfigure}\quad
  \begin{subfigure}{0.24\textwidth}
      \centering      
      \includegraphics[width=1.1\textwidth]{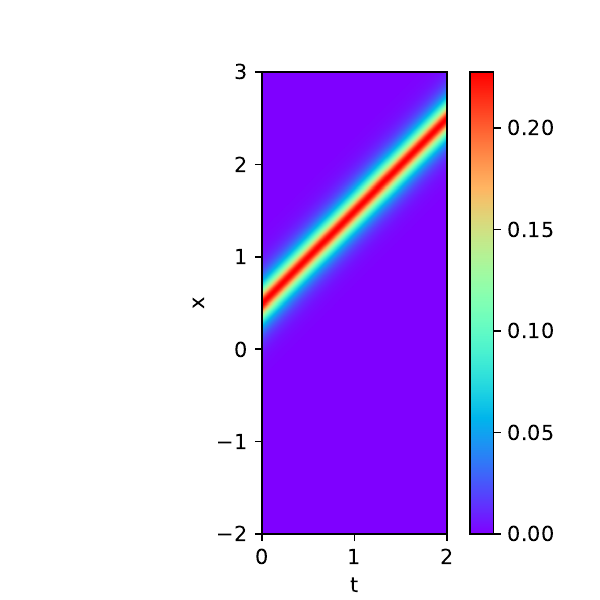}  
      \caption{Numerical solution}
      \label{fig:b}
  \end{subfigure}
  \begin{subfigure}{0.24\textwidth}
      \centering    
      \includegraphics[width=1.1\textwidth]{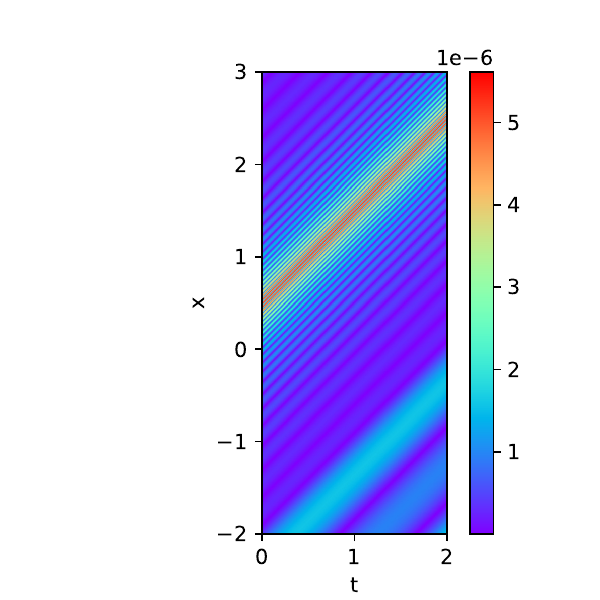}  
		 \caption{$|u_{h\tau} - u^*|$}
      \label{fig:c}
    \end{subfigure}  
  \caption{Performances of the LRNN-DG method with the wavelet function in Example \ref{ex_gKdV}.}  
  \label{figure_waveletdg}
\end{figure}

\begin{figure}[H] 
  \centering  
   
  \begin{subfigure}{0.24\textwidth}
      \centering      
      \includegraphics[width=1.1\textwidth]{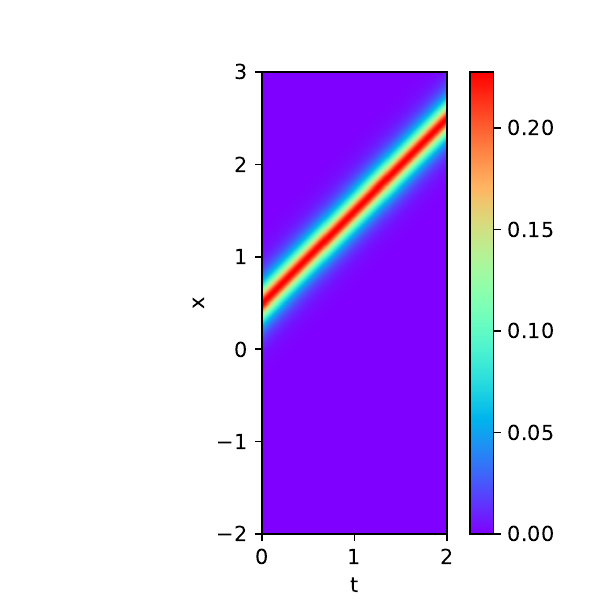}  
      \caption{Exact solution}
      \label{fig:a}
  \end{subfigure}\quad
  \begin{subfigure}{0.24\textwidth}
      \centering      
      \includegraphics[width=1.1\textwidth]{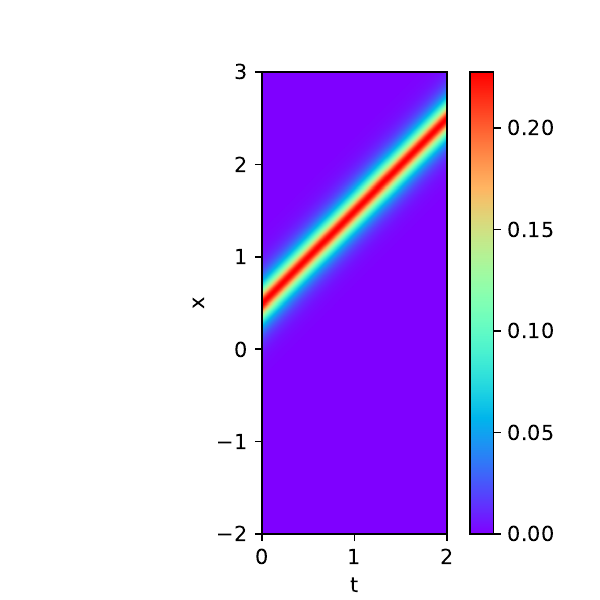}  
      \caption{Numerical solution}
      \label{fig:b}
  \end{subfigure}
  \begin{subfigure}{0.24\textwidth}
      \centering    
      \includegraphics[width=1.1\textwidth]{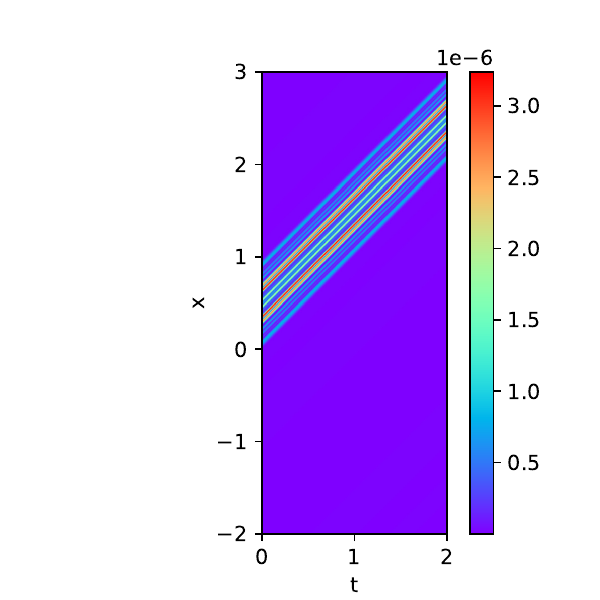}  
		 \caption{$|u_{h\tau} - u^*|$}
      \label{fig:c}
    \end{subfigure}  
  \caption{Performances of the LRNN-$C^1$DG method with the wavelet function in Example \ref{ex_gKdV}.}  
  \label{figure_waveletc1dg}
\end{figure}

After determining the slope of the characteristic direction as $k = 1.0011774546$ and the critical points $x_0 = 0.5$ of $u_0$, we define $\widehat{x} = kt + x_0 - x$ and use the Gaussian function $e^{-\widehat{x}^2 / 2}$ as the activation function. The results of the LRNN-DG methods with this activation function are provided in Table \ref{tablegKdVwavelet}, using 230 Gaussian quadrature points in each direction. For the LRNN-DG method, we set $r = 15$ and the penalty $\eta = 180/h$ or $\eta = 180/\tau$; for the LRNN-$C^1$DG method, $r = 13$ and 65 collocation points along each edge. Results for both methods with ${\rm DoF}_\sigma = 80$ are illustrated in Figures \ref{figure_waveletdg} and \ref{figure_waveletc1dg}, showing that by incorporating valuable information, the neural network with characteristic wavelet basis functions as activation functions achieves high accuracy with minimal degrees of freedom for this problem.

\begin{example}[KdV Equation with Double Solitons Collision]
\label{ex_kdv}

In this example, we examine the KdV equation that models the collision of two solitons:
\[
u_t + u u_x + \epsilon u_{xxx} = 0 \quad (t, x) \in I \times \Omega,
\]
with the initial condition defined as:
\[
u_0(x) = 3c_1 \,  {\rm sech}^2 (k_1(x - x_1)) + 3c_2 \,  {\rm sech}^2 (k_2(x - x_2)),
\]
where $c_1 = 0.3$, $c_2 = 0.1$, $x_1 = 0.4$, $x_2 = 0.8$, and $k_i = \frac{1}{2} \sqrt{\frac{c_i}{\epsilon}}$ for $i = 1, 2$. Here, $\epsilon = 4.84 \times 10^{-4}$. The spatial domain is $\Omega = (0, 2)$ and the time interval is $I = (0, 2)$, with periodic boundary conditions applied. While no analytical solution exists for this problem, a numerical solution computed using the DG method is available in \cite{Yan2002dgkdv}.
\end{example}

\begin{figure}[H] 
  \centering  
   
  \begin{subfigure}{0.205\textwidth}
      \centering      
      \includegraphics[width=1.15\textwidth]{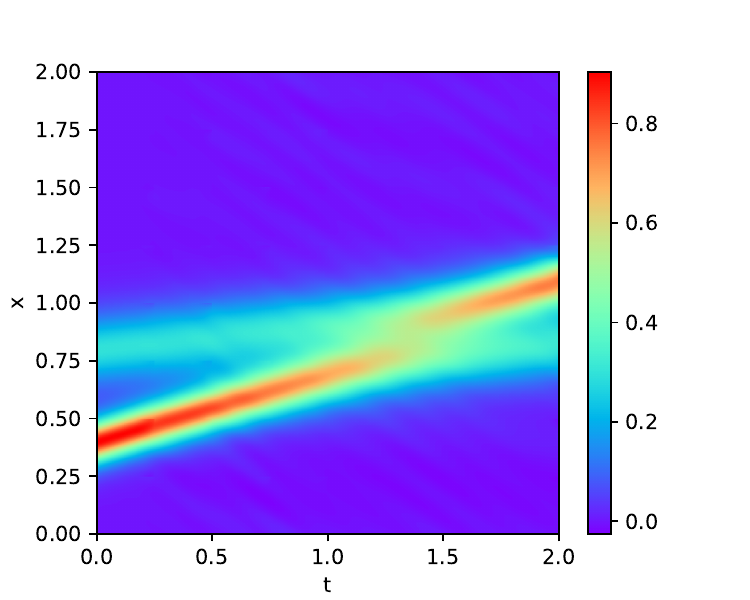}  
      \caption{$u_{h\tau}$ with 64 elements}
      %\label{fig:a}
  \end{subfigure}
  \begin{subfigure}{0.28\textwidth}
      \centering      
      \includegraphics[width=1\textwidth]{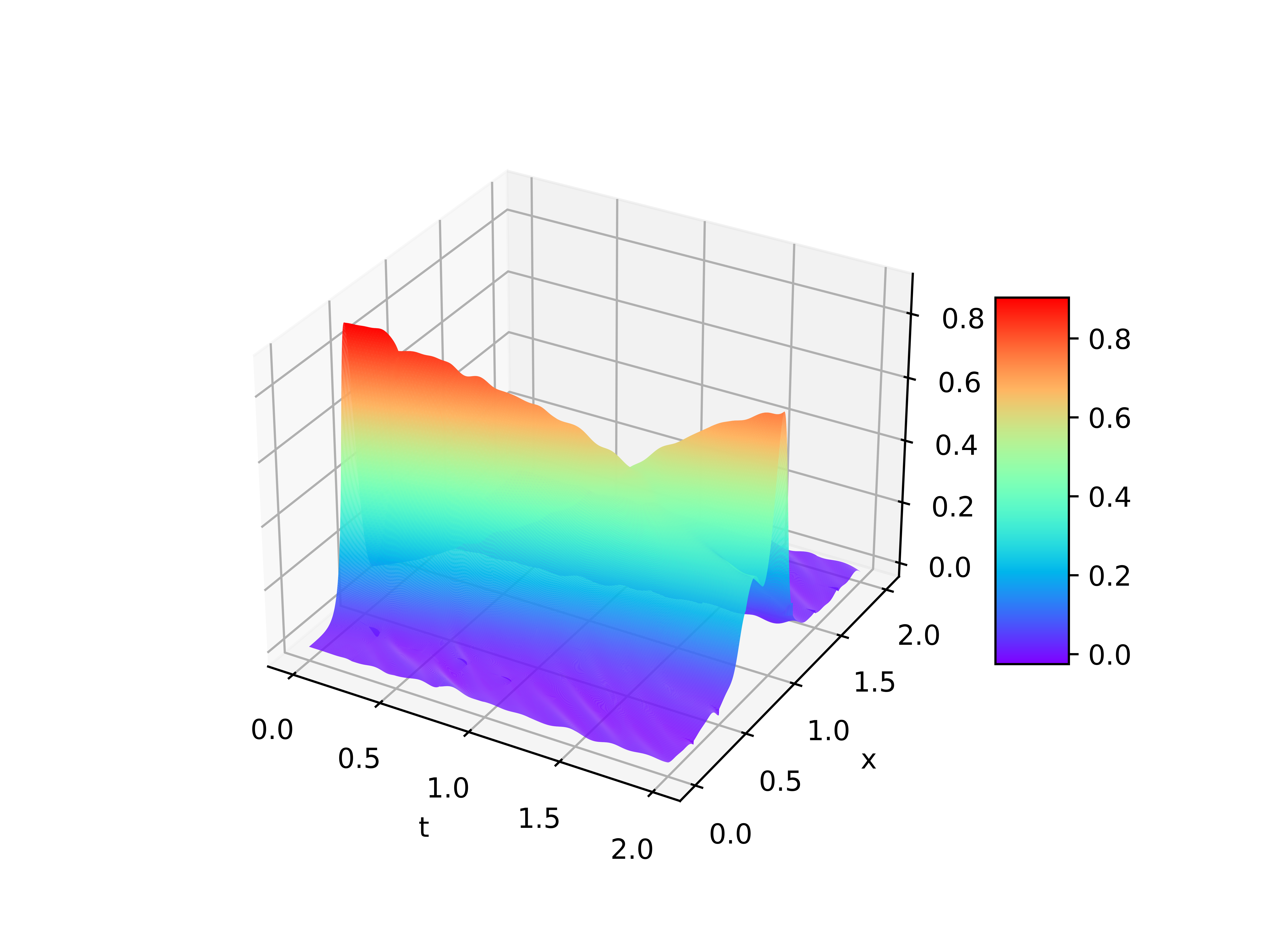}  
      \caption{$u_{h\tau}$ with 64 elements}
      %\label{fig:c}
  \end{subfigure}\;
  \begin{subfigure}{0.205\textwidth}
      \centering    
      \includegraphics[width=1.15\textwidth]{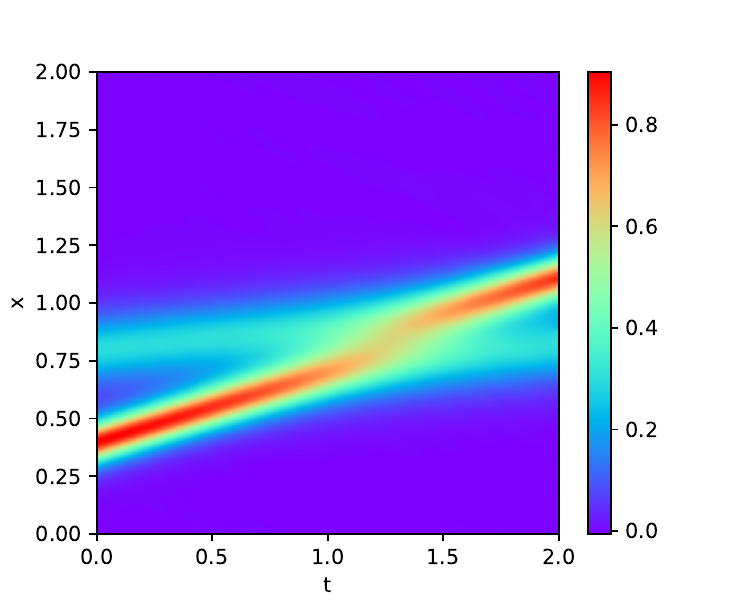}  
		 \caption{$u_{h\tau}$ with 81 elements}
      %\label{fig:b}
    \end{subfigure}  
   \begin{subfigure}{0.28\textwidth}
      \centering    
      \includegraphics[width=1\textwidth]{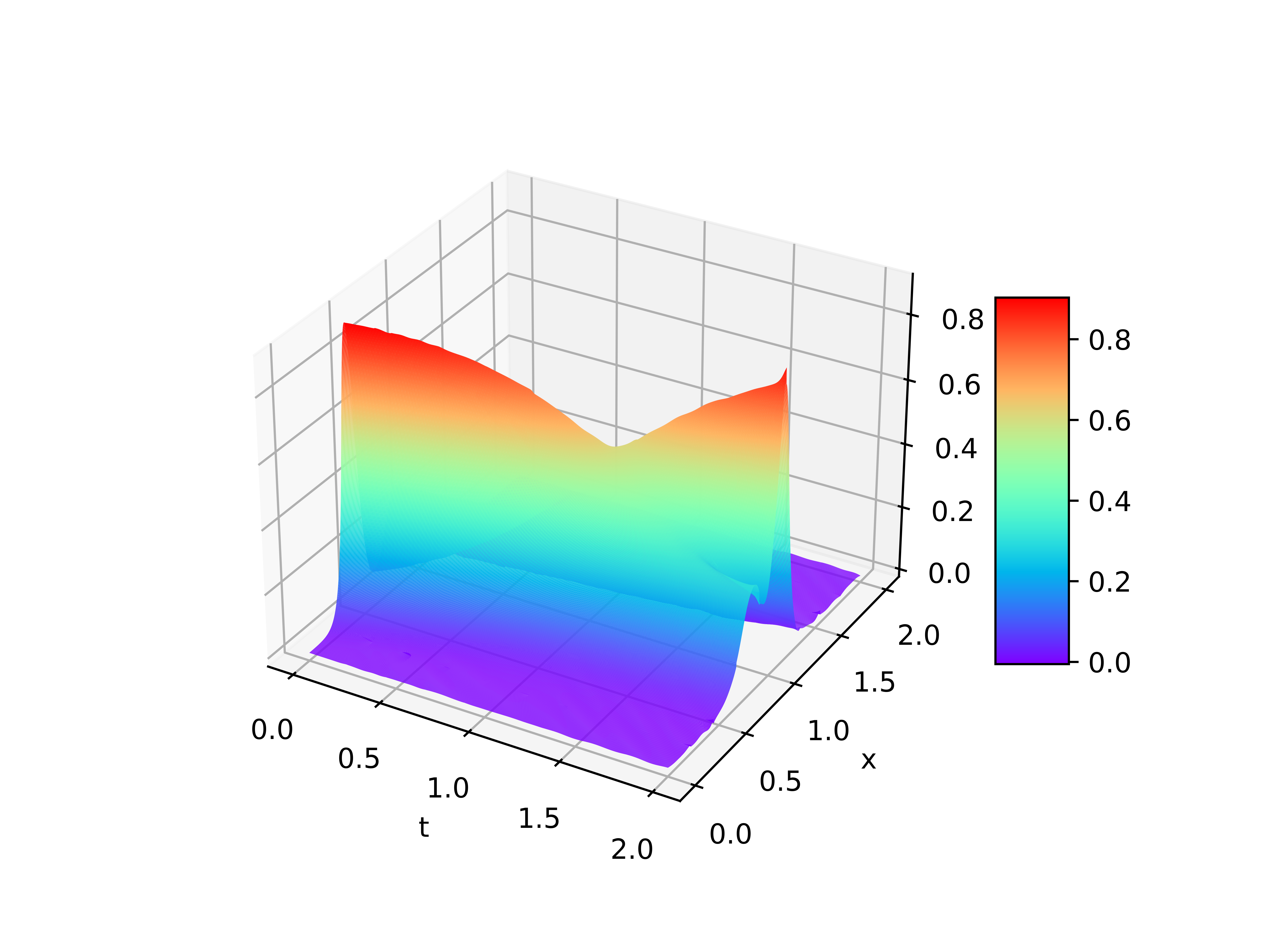}  
		 \caption{$u_{h\tau}$ with 81 elements}
      %\label{fig:d}
    \end{subfigure}  
  \caption{ Numerical solutions of the LRNN-$C^1$DG method on uniform meshes in Example \ref{ex_kdv}.}  
  \label{kdv_c1dg_dsc}
\end{figure}

In this experiment, we first apply the LRNN-$C^1$DG method on a uniform mesh and display the space-time solution with mesh sizes $8 \times 8$ and $9 \times 9$, as shown in Figure \ref{kdv_c1dg_dsc}. The parameters are set to $r = 1$, $\theta = 0.4$, $N_{ni} = 30$, $\epsilon_0 = 10^{-4}$, with 17 collocation points on each edge and 10 Gaussian quadrature points in each direction. Despite the dense refinement of the mesh, the numerical solution still fails to capture some of the finer details accurately.

\vspace{-6mm}

\begin{figure}[H] 
  \centering  
   
  \begin{subfigure}{0.3\textwidth}
      \centering      
      \includegraphics[width=0.88\textwidth]{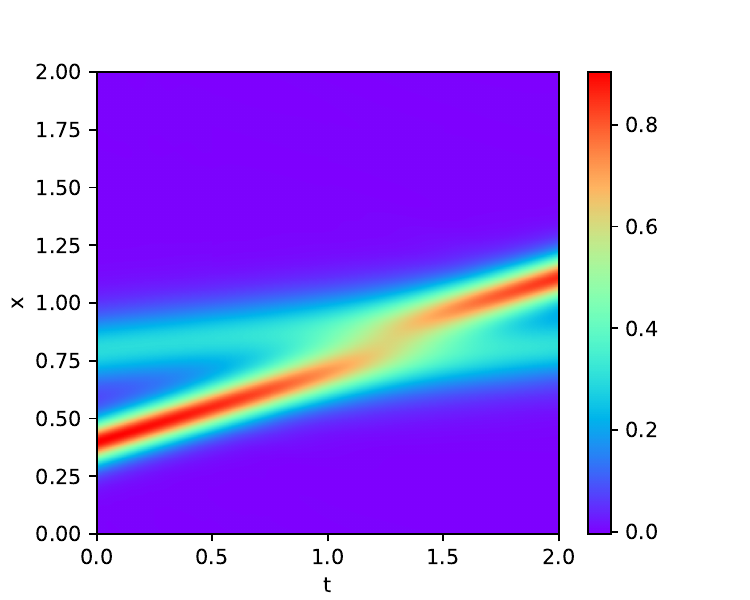}  
      \caption{$u_{h\tau}$ with 43 elements}
      %\label{fig:a}
  \end{subfigure}
  \begin{subfigure}{0.33\textwidth}
      \centering      
      \includegraphics[width=0.92\textwidth]{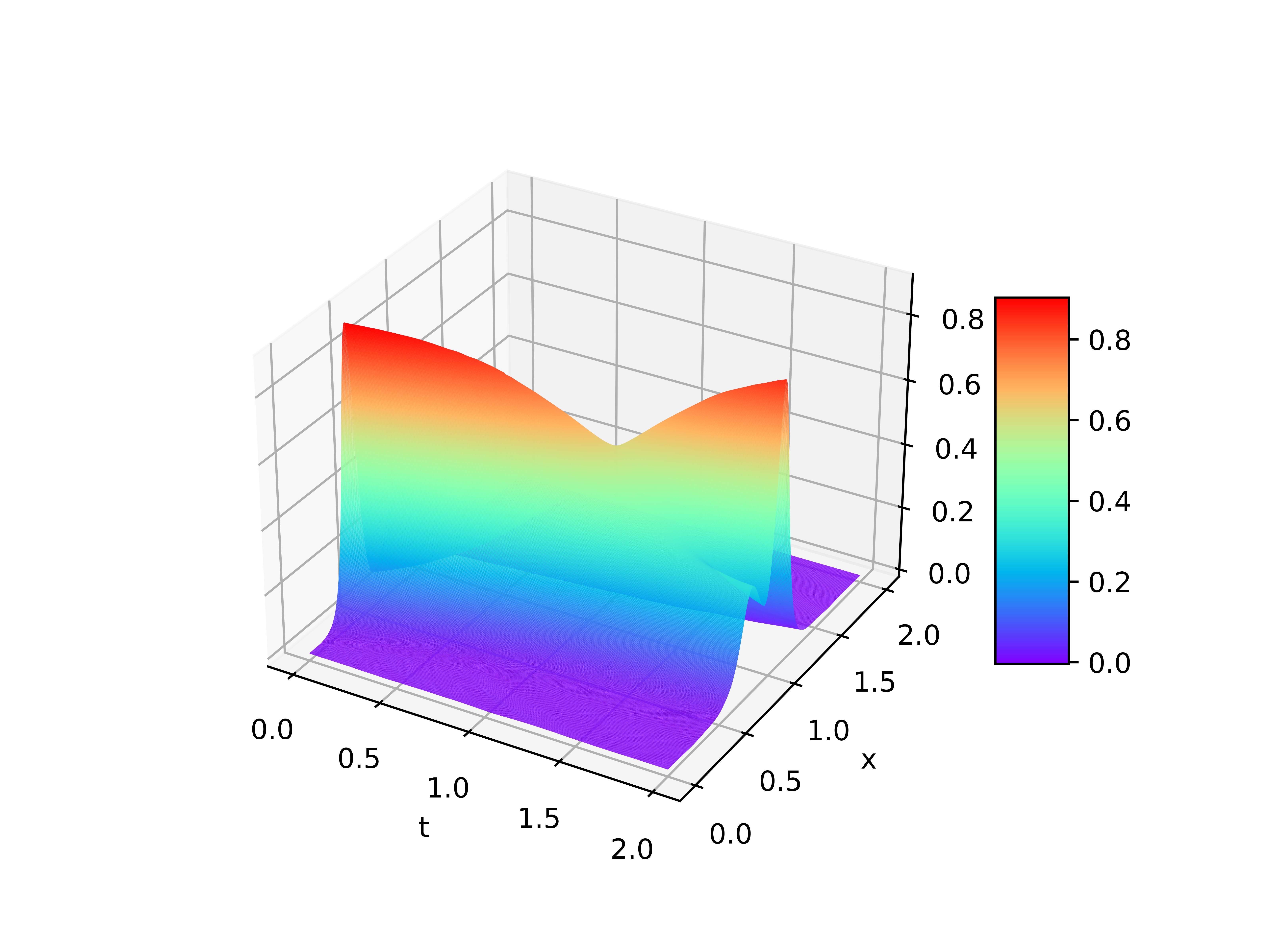}  
      \caption{$u_{h\tau}$ with 43 elements}
      %\label{fig:c}
  \end{subfigure}
  \begin{subfigure}{0.3\textwidth}
      \centering      
      \includegraphics[width=0.88\textwidth]{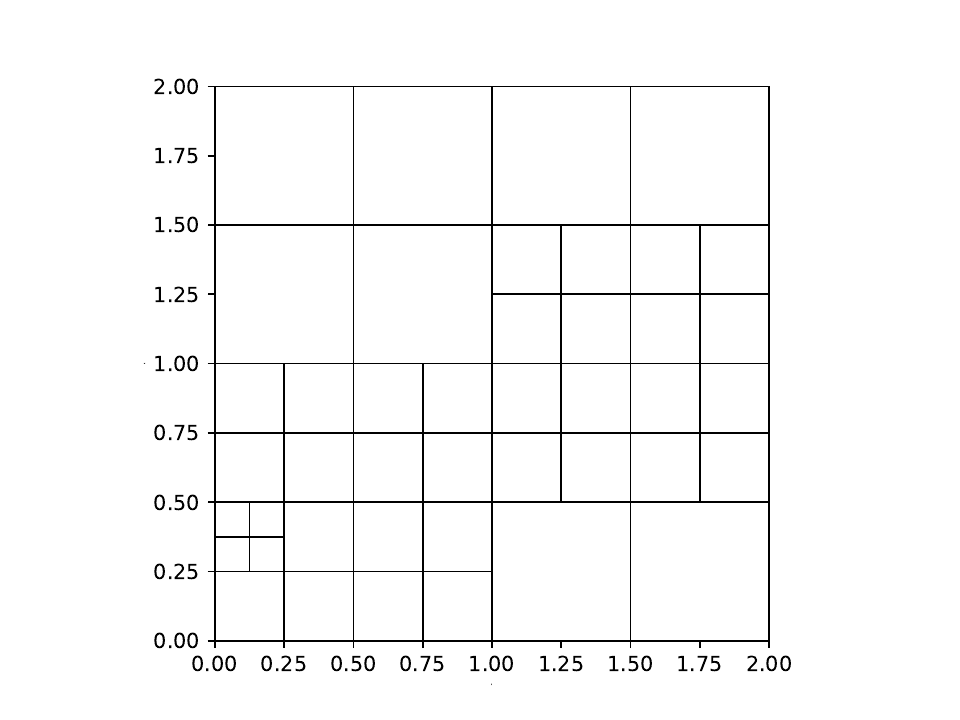}  
      \caption{Adaptive mesh with 43 elements}
      \label{fig:e}
  \end{subfigure}\\
  \begin{subfigure}{0.3\textwidth}
      \centering    
      \includegraphics[width=0.88\textwidth]{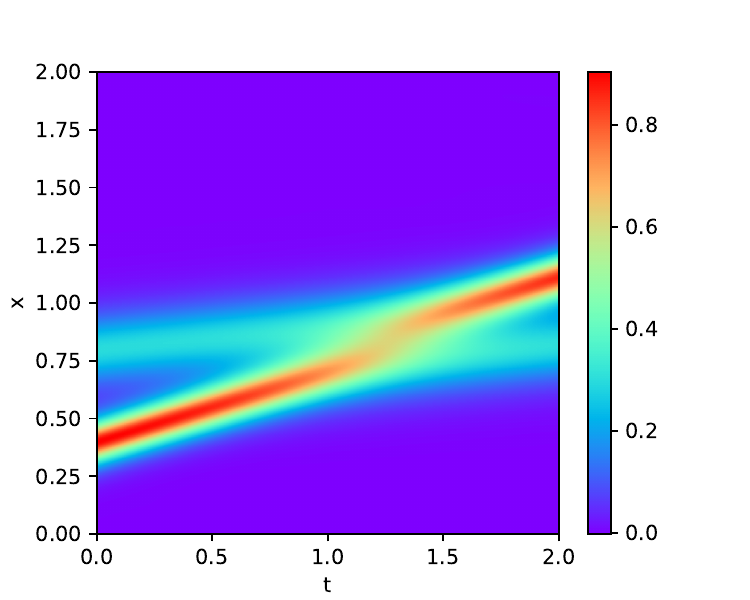}  
		 \caption{$u_{h\tau}$ with 58 elements}
      \label{fig:b}
    \end{subfigure} 
   \begin{subfigure}{0.33\textwidth}
      \centering    
      \includegraphics[width=0.92\textwidth]{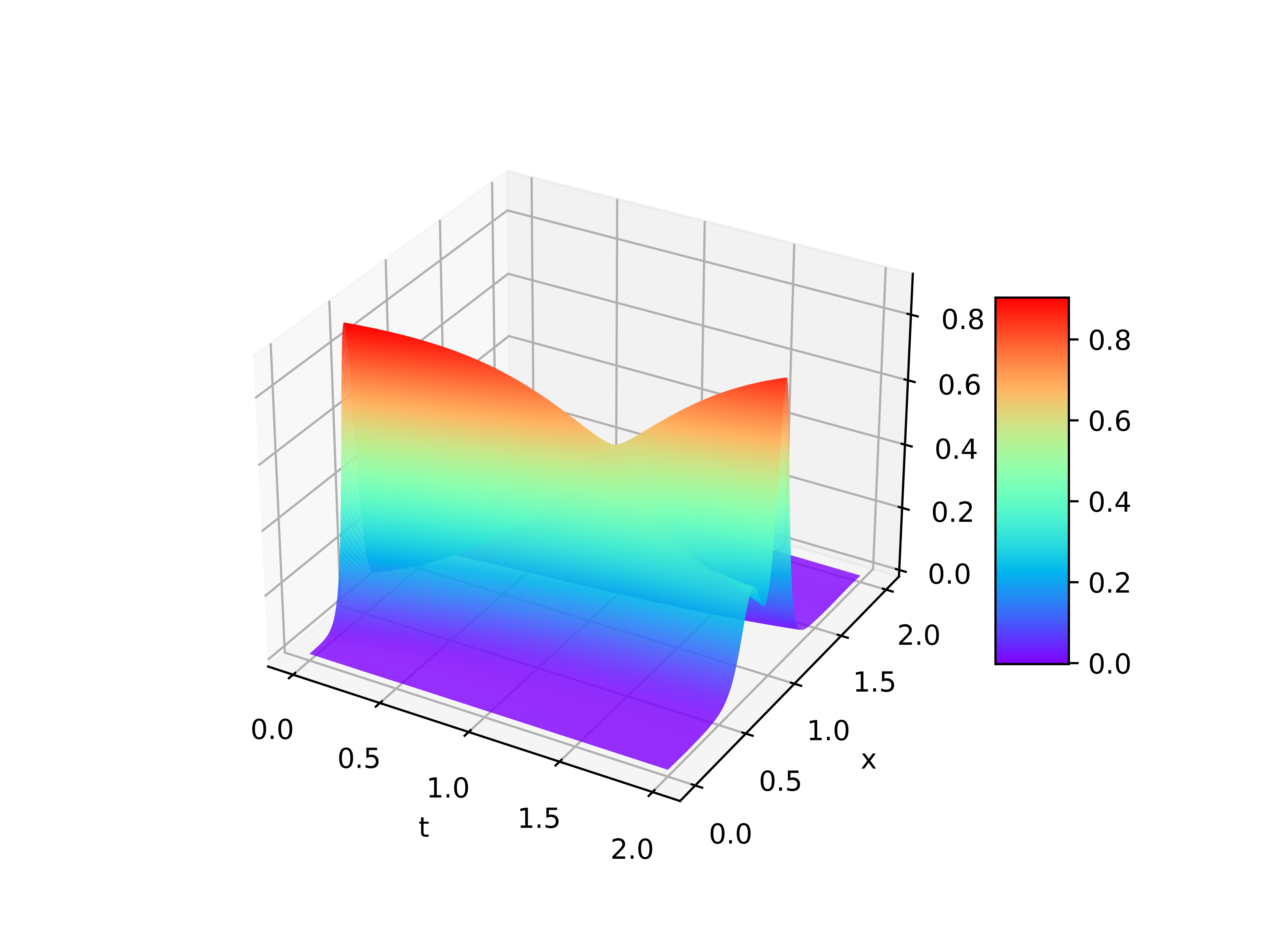}  
		 \caption{$u_{h\tau}$ with 58 elements}
      \label{fig:d}
    \end{subfigure} 
   \begin{subfigure}{0.3\textwidth}
      \centering    
      \includegraphics[width=0.88\textwidth]{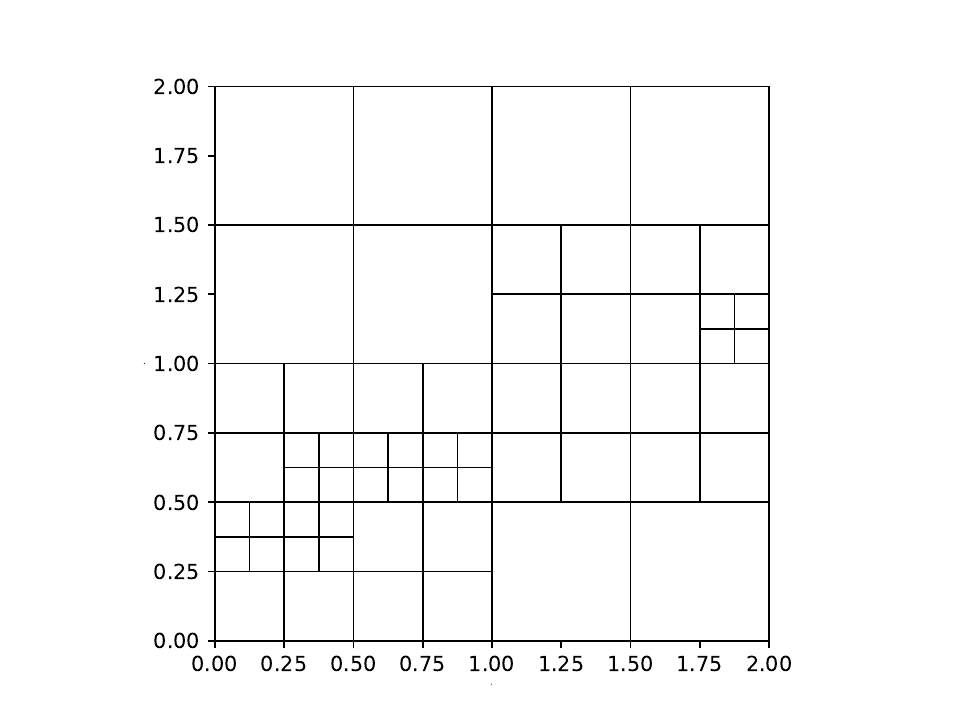}  
		 \caption{Adaptive mesh with 58 elements}
      \label{fig:f}
    \end{subfigure}  
  \caption{Numerical solutions of the LRNN-$C^1$DG method on adaptive meshes in Example \ref{ex_kdv}.}  
  \label{kdv_adac1dg_dsc}
\end{figure}

Subsequently, we adopt the same approach on an adaptive mesh to tackle the challenge of the collision of double solitons. Starting with a $4\times 4$ uniform initial mesh, we retain the other settings used in the uniform mesh approach. The outcomes of two adaptive refinements result in 43 and 58 elements, respectively, as shown in Figure \ref{kdv_adac1dg_dsc}. This demonstrates that the LRNN-$C^1$DG method, when applied with an adaptive mesh approach, can effectively capture the details of the solution.

\begin{example}[2D Burgers Equation]
\label{ex_bur}

In this example, we investigate the 2D Burgers equation:
\begin{align*}
u_t + u(u_x + u_y) - \epsilon \Delta u &=0 \quad (t,x,y)\in I\times\Omega,\\
u(t,x,y) &= g(t,x,y)\quad (t,x,y)\in I \times {\Gamma_D},\\
u(0, x, y) &= u_0(x,y) \quad (x,y)\in \Omega,
\end{align*}
where the exact solution is given by $u = 1/(1+e^{\frac{x+y-t}{2\epsilon}})$.
The spatial domain is $\Omega =(0,1)$, the time interval is $I=(0,1)$, and the boundary $\Gamma_D  = \partial \Omega$. The boundary function $g$ and initial condition $u_0$ are derived from the exact solution.
\end{example}

\begin{table}[h]
\centering
\begin{tabular}{|c|c|l|l|l|}
\hline
                              &                                              & \multicolumn{1}{c|}{}                     & \multicolumn{1}{c|}{}                      & \multicolumn{1}{c|}{}                      \\
\multirow{-2}{*}{$\tau$, $h$} & \multirow{-2}{*}{\diagbox[width=5em,trim=l]{Norm}{${\rm DoF}_\sigma$}} & \multicolumn{1}{c|}{\multirow{-2}{*}{80}} & \multicolumn{1}{c|}{\multirow{-2}{*}{160}} & \multicolumn{1}{c|}{\multirow{-2}{*}{320}} \\ \hline
                              & $E^{L^{2}}$                                        & 1.20E-02                                  & 5.78E-05                                   & 8.58E-06                                   \\ \cline{2-5} 
                              & $E^{H^{1}}$                                        & 3.28E-01                                  & 2.36E-03                                   & 4.00E-04                                   \\ \cline{2-5} 
\multirow{-3}{*}{1/2}         & $E^{L^{2}}(t=1)$                                  & 1.72E-02                                  & 7.34E-05                                   & 1.37E-05                                   \\ \hline
                              & $E^{L^{2}}$                                        & 3.30E-03                                  & 8.12E-06                                   & 8.69E-07                                   \\ \cline{2-5} 
                              & $E^{H^{1}}$                                        & 1.37E-01                                  & 4.82E-04                                   & 6.07E-05                                   \\ \cline{2-5} 
\multirow{-3}{*}{1/3}         & $E^{L^{2}}(t=1)$                                  & 2.94E-03                                  & 1.21E-05                                   & 1.28E-06                                   \\ \hline
                              & $E^{L^{2}}$                                        & 4.84E-04                                  & 3.56E-06                                   & 4.22E-07                                   \\ \cline{2-5} 
                              & $E^{H^{1}}$                                        & 2.63E-02                                  & 2.91E-04                                   & 3.89E-05                                   \\ \cline{2-5} 
\multirow{-3}{*}{1/4}         & $E^{L^{2}}(t=1)$                                  & 7.21E-04                                  & 4.95E-06                                   & 5.18E-07                                   \\ \hline
\end{tabular}
\caption{Errors of the space-time LRNN-DG method in Example \ref{ex_bur} when $\epsilon = 0.1$}
\label{tablegbur01dg}
\end{table}

\begin{table}[h]
\centering
\begin{tabular}{|c|c|l|l|l|}
\hline
$\tau$, $h$          & \diagbox[width=5em,trim=l]{Norm}{${\rm DoF}_\sigma$} & \multicolumn{1}{c|}{80} & \multicolumn{1}{c|}{160} & \multicolumn{1}{c|}{320} \\ \hline
\multirow{3}{*}{1/2} & $E^{L^{2}}$                    & 8.68E-03                & 3.20E-04                 & 2.02E-04                 \\ \cline{2-5} 
                     & $E^{H^{1}}$                    & 1.06E-01                & 8.49E-03                 & 3.89E-03                 \\ \cline{2-5} 
                     & $E^{L^{2}} (t=1)$              & 1.57E-02                & 6.46E-04                 & 2.92E-04                 \\ \hline
\multirow{3}{*}{1/3} & $E^{L^{2}}$                    & 2.41E-03                & 5.04E-05                 & 7.77E-06                 \\ \cline{2-5} 
                     & $E^{H^{1}}$                    & 4.30E-02                & 1.52E-03                 & 2.06E-04                 \\ \cline{2-5} 
                     & $E^{L^{2}} (t=1)$              & 4.91E-03                & 1.28E-04                 & 2.60E-05                 \\ \hline
\multirow{3}{*}{1/4} & $E^{L^{2}}$                    & 7.62E-04                & 8.77E-06                 & 1.03E-06                 \\ \cline{2-5} 
                     & $E^{H^{1}}$                    & 1.65E-02                & 4.74E-04                 & 3.32E-05                 \\ \cline{2-5} 
                     & $E^{L^{2}} (t=1)$              & 1.43E-03                & 6.33E-05                 & 2.59E-06                 \\ \hline
\end{tabular}
\caption{Errors of the space-time LRNN-$C^1$DG method in Example \ref{ex_bur} when $\epsilon = 0.1$}
\label{tablegbur01c1dg}
\end{table}

\begin{table}[h]
\centering
\begin{tabular}{|c|c|l|l|}
\hline
                              &                                              & \multicolumn{1}{c|}{}                      & \multicolumn{1}{c|}{}                      \\
\multirow{-2}{*}{$\tau$, $h$} & \multirow{-2}{*}{\diagbox[width=5em,trim=l]{Norm}{${\rm DoF}_\sigma$}} & \multicolumn{1}{c|}{\multirow{-2}{*}{160}} & \multicolumn{1}{c|}{\multirow{-2}{*}{320}} \\ \hline
                              & $E^{L^{2}}$                                        & 3.89E-02                                   & 2.30E-02                                   \\ \cline{2-4} 
\multirow{-2}{*}{1/3}         & $E^{L^{2}} (t=1)$                                  & 5.08E-02                                   & 3.58E-02                                   \\ \hline
                              & $E^{L^{2}}$                                        & 2.32E-02                                   & 1.27E-02                                   \\ \cline{2-4} 
\multirow{-2}{*}{1/4}         & $E^{L^{2}} (t=1)$                                  & 3.37E-02                                   & 2.03E-02                                   \\ \hline
                              & $E^{L^{2}}$                                        & 1.45E-02                                   & 7.46E-03                                   \\ \cline{2-4} 
\multirow{-2}{*}{1/5}         & $E^{L^{2}} (t=1)$                                  & 2.16E-02                                   & 1.19E-02                                   \\ \hline
                              & $E^{L^{2}}$                                        & 9.47E-03                                   & 4.19E-03                                   \\ \cline{2-4} 
\multirow{-2}{*}{1/6}         & $E^{L^{2}} (t=1)$                                  & 1.30E-02                                   & 7.37E-03                                   \\ \hline
\end{tabular}
\caption{Errors of the space-time LRNN-DG method in Example \ref{ex_bur} when $\epsilon = 0.01$}
\label{tablegbur001dg}
\end{table}

\begin{table}[h]
\centering
\begin{tabular}{ccccccc}
\toprule
{$N_e$}  & \multicolumn{1}{c}{8} & \multicolumn{1}{c}{22} & \multicolumn{1}{c}{36} & \multicolumn{1}{c}{85} & \multicolumn{1}{c}{120} & 148                          \\ \midrule
$E^{L^2}$      & 6.10E-02              & 4.73E-02               & 2.30E-02               & 1.67E-02               & 1.21E-02                & \multicolumn{1}{l}{5.58E-03} \\
$E^{L^2}(t=1)$ & 1.07E-01              & 5.87E-02               & 3.11E-02               & 2.54E-02               & 2.19E-02                & 8.99E-03           \\ \bottomrule
\end{tabular}
\caption{Errors of the adaptive LRNN-DG method in Example \ref{ex_bur} when $\epsilon = 0.01$}
\label{tablegbur001adadg}
\end{table}

To compare LRNN-DG methods with the RNN-PG approach from \cite{Shang2022DeepPetrov}, we report the $L^2$ errors of the proposed method at $t = 1$. In Table \ref{tablegbur01dg}, the errors of the LRNN-DG method for $\epsilon = 0.1$ are presented, where the parameters are set to $r = 0.6$ and the interior penalty $\eta = 40/h$ or $\eta = 40/\tau$. Similarly, Table \ref{tablegbur01c1dg} shows the errors for the LRNN-$C^1$DG method with $\epsilon = 0.1$, using parameters $r = 0.4$ and 20 collocation points along each edge. For both methods, $N_{ni} = 25$, $\epsilon_0 = 10^{-4}$, and 12 Gaussian integration points are used in each direction. Notably, the LRNN-DG method outperforms the LRNN-$C^1$DG method, with both demonstrating superior accuracy compared to the RNN-PG method from \cite{Shang2022DeepPetrov}.

For $\epsilon = 0.01$, achieving satisfactory results with the RNN-PG method is challenging. To address this, we apply the LRNN-DG method on a uniform mesh to solve the problem for $\epsilon = 0.01$, as detailed in Table \ref{tablegbur001dg}. Parameters include $r = 0.7$, interior penalty $\eta = 45/h$ or $\eta = 45/\tau$, 12 Gaussian integration points per direction, along with $N_{ni} = 25$ and $\epsilon_0 = 10^{-4}$.

\begin{figure}[H] 
  \centering  
  \begin{subfigure}{0.3\textwidth}
      \centering      
      \includegraphics[width=0.9\textwidth]{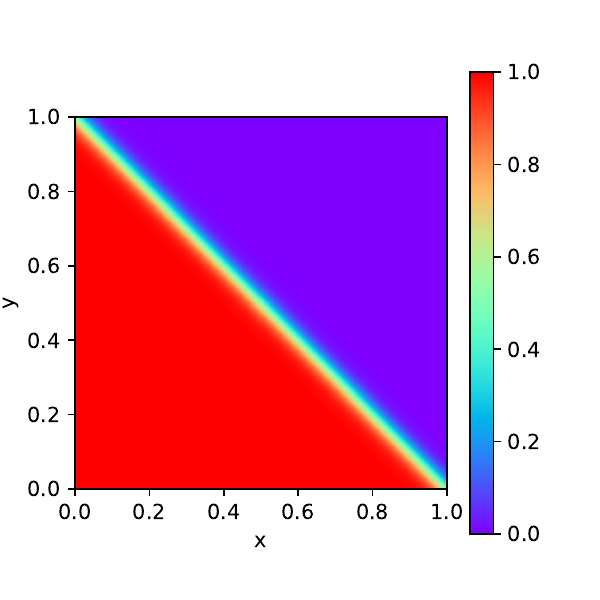}  
      \caption{Exact solution}
      %\label{fig:a}
  \end{subfigure}\quad
  \begin{subfigure}{0.3\textwidth}
      \centering      
      \includegraphics[width=0.9\textwidth]{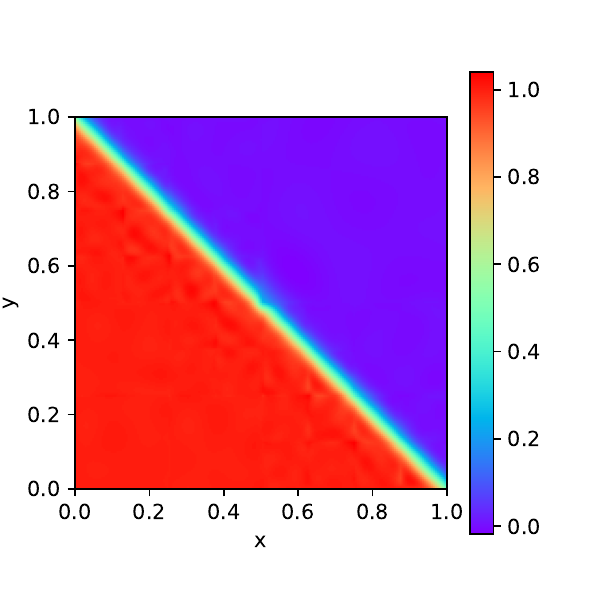}  
      \caption{Numerical solution}
      %\label{fig:b}
  \end{subfigure}\quad
  \begin{subfigure}{0.3\textwidth}
      \centering    
      \includegraphics[width=0.9\textwidth]{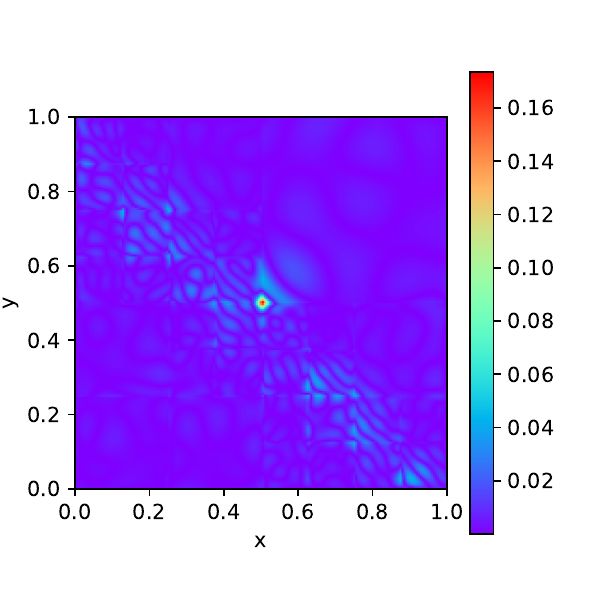}  
		 \caption{Absolute errors}
      %\label{fig:c}
    \end{subfigure}  \\
   \begin{subfigure}{0.4\textwidth}
      \centering    
      \includegraphics[width=1\textwidth]{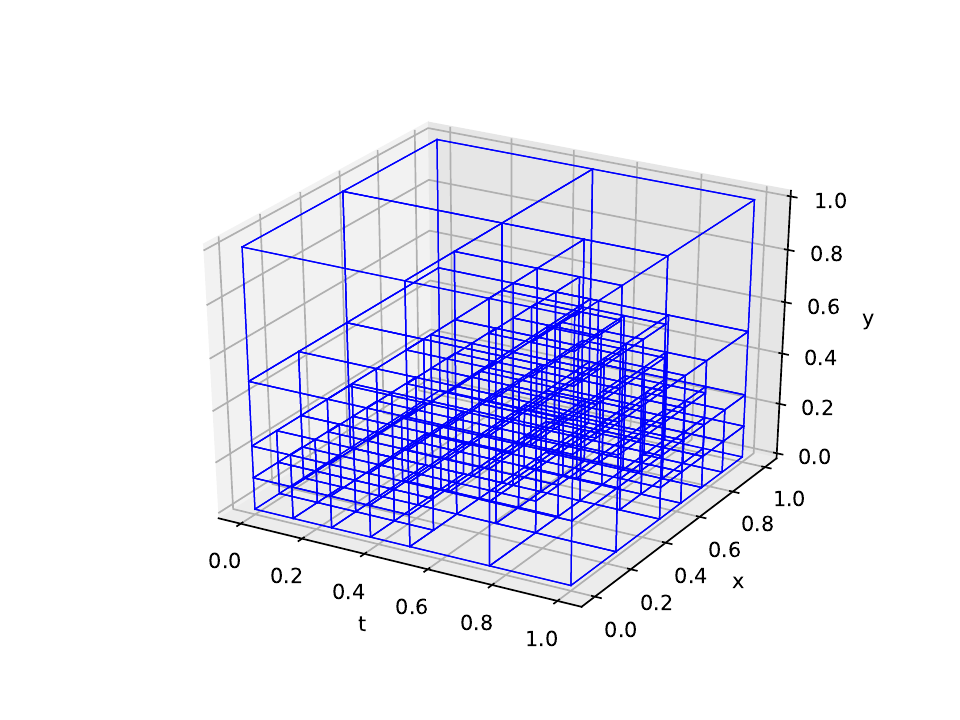}  
		 \caption{Adaptive mesh of LRNN-DG method}
      %\label{fig:d}
    \end{subfigure}  \quad\quad\quad
    \begin{subfigure}{0.4\textwidth}
      \centering    
      \includegraphics[width=1\textwidth]{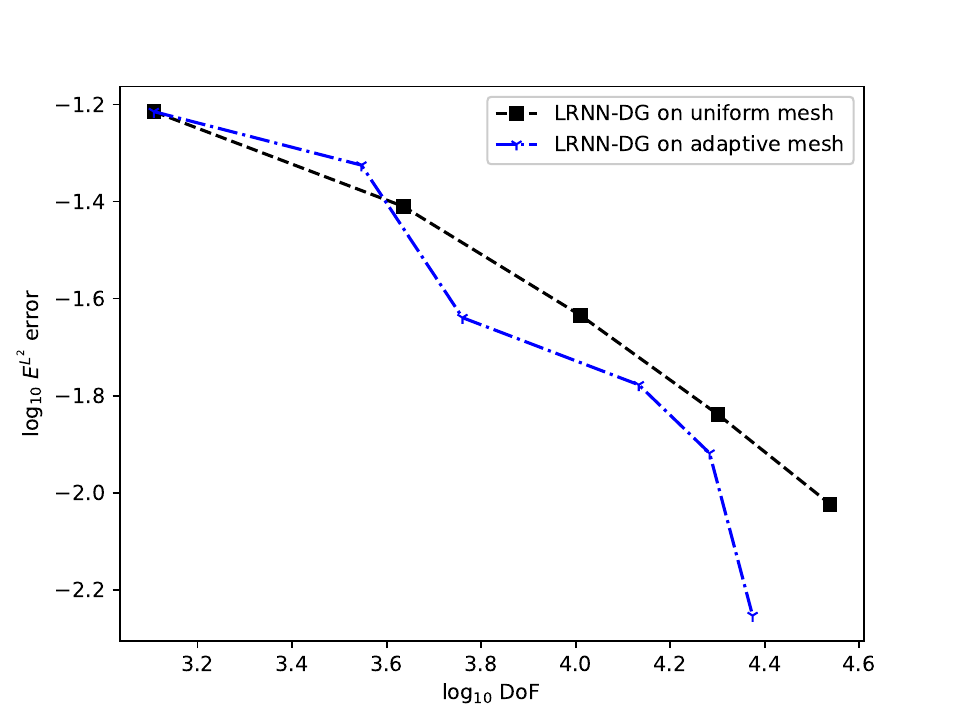}  
		 \caption{$E^{L^2}$ errors with respect to DoF}
      \label{fig:bur_comparison}
    \end{subfigure}  
  \caption{The performances of the adaptive LRNN-DG method when $\epsilon = 0.01$ in Example \ref{ex_bur}.}  
  \label{figure_adapdg_bur}
\end{figure}

Additionally, we apply the LRNN-DG method on an adaptive mesh for the same scenario with $\epsilon = 0.01$, as detailed in Table \ref{tablegbur001adadg}. Using the same parameter settings and ${\rm DoF}_\sigma = 160$, the adaptive approach is implemented. Figure \ref{figure_adapdg_bur} illustrates the exact solution, numerical solution, absolute errors at $t=1$, and the adaptive mesh decomposition for the LRNN-DG method with ${\rm DoF}_\sigma = 160$. Finally, we compare the performance of the LRNN-DG method on both uniform and adaptive meshes in Figure \ref{fig:bur_comparison}, which highlights the $E^{L^2}$ error when ${\rm DoF}_\sigma = 160$, with all other parameters consistent with previous tables.

\section{Summary}
\label{summary}

In this study, we employed space-time LRNN-DG methods to solve nonlinear equations, including KdV-type equations involving single or double soliton collisions and Burgers' equations with a small coefficient for the second derivative term. To enhance the approximation capabilities of LRNN-DG methods, we explored two alternative strategies: adaptive mesh refinement, which uses error indicators to guide mesh adaptation, and characteristic meshes, which incorporate prior information to reduce solution complexity within each subdomain. Our experiments demonstrated the effectiveness of LRNN-DG methods on both adaptive and characteristic meshes, revealing several key advantages:  
(i) LRNN-DG methods can achieve highly accurate numerical solutions with relatively few degrees of freedom;  
(ii) LRNN-DG methods on adaptive or characteristic meshes require fewer degrees of freedom to attain higher accuracy compared to uniform meshes;  
(iii) Space-time LRNN-DG methods efficiently solve time-dependent problems, mitigating error accumulation over time.

While LRNN-DG methods show great promise, they also pose challenges and open questions for future research. A thorough numerical analysis of these methods is essential to establish their theoretical foundations. Furthermore, developing reliable and efficient error estimators is a compelling direction for future work. Lastly, exploring the construction of characteristic meshes to further enhance computational efficiency remains a promising avenue for future studies.

\end{document}